\documentclass[journal]{IEEEtran}

\ifCLASSOPTIONcompsoc
\else
\fi
\ifCLASSINFOpdf
\else
\fi
\usepackage[cmex10]{amsmath}
\usepackage{amssymb}
\usepackage{enumerate}
\usepackage{array}
\usepackage{fixltx2e}

\usepackage{stfloats}
\newcommand\MYhyperrefoptions{bookmarks=true,bookmarksnumbered=true,
pdfpagemode={UseOutlines},plainpages=false,pdfpagelabels=true,
colorlinks=true,linkcolor={black},citecolor={black},pagecolor={black},
urlcolor={black},
pdftitle={A Viterbi process for general hidden Markov models},
pdfsubject={Typesetting},
pdfauthor={J\"uri  Lember},
pdfkeywords={asymptotic, HMM, MAP sequence, Viterbi algorithm, Viterbi training, Viterbi extraction}}
\ifCLASSINFOpdf
\usepackage[\MYhyperrefoptions,pdftex]{hyperref}
\else
\usepackage[\MYhyperrefoptions,breaklinks=true,dvips]{hyperref}
\usepackage{breakurl}
\fi

\newtheorem{theorem}{Theorem}[section]
\newtheorem{definition}[theorem]{Definition}
\newtheorem{example}[theorem]{Example}
\def\V{\mathcal{V}}
\def\P{{\mathbf P}}
\newtheorem{lemma}{Lemma}[section]
\def\W{\mathcal{W}}
\begin{document}
%
\title{A constructive proof of the existence of Viterbi processes}
%
%
%

\author{J\"uri~Lember,~
        Alexey~Koloydenko~
\IEEEcompsocitemizethanks{
\IEEEcompsocthanksitem J. Lember is with the Institute 
of Mathematical Statistics, Tartu University, J. Liivi 2-507, 
50409, Estonia.\protect\\
E-mail: jyri@ut.ee
\IEEEcompsocthanksitem A. Koloydenko is with the Division of Statistics of Nottingham University, 
University Park, Nottingham, NG7 2RD, UK.\protect\\
E-mail: alexey.koloydenko@nottingham.ac.uk}
\thanks{Manuscript received April 8, 2008; revised ??, 2008
}}

\markboth{Submitted to the IEEE Transactions on Information Theory }
{Lember \MakeLowercase{\textit{et al.}}: The Viterbi Process: Construction and Properties}
%

\IEEEcompsoctitleabstractindextext{%
\begin{abstract}
Since the early days of digital communication, hidden Markov models
(HMMs) have now been also routinely used in speech recognition, processing of
natural languages, images, and in bioinformatics.  In an HMM
$(X_i,Y_i)_{i\ge 1}$,  observations $X_1,X_2,\ldots$ are assumed to be
conditionally independent given an ``explanatory'' Markov process
$Y_1,Y_2,\ldots$, which itself is not observed; moreover, the
conditional distribution of $X_i$ depends solely on $Y_i$. Central
to the theory and applications of HMM is the Viterbi algorithm to find
{\em a maximum a posteriori} (MAP) estimate $q_{1:n}=(q_1,q_2,\ldots,q_n)$ of
$Y_{1:n}$ given  observed data $x_{1:n}$.
Maximum {\em a posteriori} paths are also known as
Viterbi paths or  alignments.  Recently, attempts have been made to study the behavior of
Viterbi alignments when $n\to \infty$. Thus, it has  been
shown that in some special cases a well-defined limiting Viterbi alignment exists. While
innovative, these attempts have relied on rather strong assumptions 
and involved proofs which are existential.  This work proves the
existence of infinite Viterbi alignments in a more constructive manner and 
for a very general class of HMMs.
\end{abstract}

}

\maketitle

\IEEEdisplaynotcompsoctitleabstractindextext

%
\IEEEpeerreviewmaketitle

\section{Introduction}\label{sec:intro}
\IEEEPARstart{L}{et} $Y=(Y_i)_{i\ge 1}$ be a Markov chain  with 
state space $S=\{1,\ldots,K\}$, $K>1$, and transition matrix $\mathbb{P}=(p_{ij})_{i,j\in S}$. 
Suppose that $Y$ is irreducible and aperiodic, hence a
unique stationary distribution $\pi=\pi\mathbb{P}$ exists; 
suppose further that $Y_i\sim \pi$ from time $i=1$.
To every state $l\in S$, let us assign an
{\it emission distribution} $P_l$
on $({\cal X},{\cal B})$, where 
${\cal X}=\mathbb{R}^D$, the $D$-dimensional Euclidean space.
Let $f_l$ be the density of $P_l$ with respect to a suitable
reference measure $\lambda$ on $({\cal X}, {\cal B})$. Most commonly,
$\lambda$ is either the Lebesgue measure
(continuously distributed $X_i$) or the counting measure
(discretely distributed $X_i$).

\begin{definition}\label{def:HMM} The stochastic process $(X,Y)$ is a
hidden Markov model if there is a (measurable) function $h$ such
that for each $n$,
$X_n=h(Y_n,e_n)$, where $e_1,e_2,\ldots$
are i.i.d. and independent of $Y$.
\end{definition}
Hence, the emission distribution $P_l$ is the distribution of
$h(l,e_n)$. The distribution of $X$ is completely determined by 
$\mathbb{P}$ and the emission distributions
$P_l,$ $l\in S$. It can be shown that $X$ is also ergodic \cite{HMP, volatility, leroux}.
Let $x_{1:n}=(x_1,\ldots,x_n)$ and $y_{1:n}=(y_1,\ldots,y_n)$
be fixed observed and unobserved realizations, respectively, 
of HMM $(X_i,Y_i)_{i\ge 1}$ up to time $n$. 
Treating $y_{1:n}$ as parameters to be estimated, let
$\Lambda(q_{1:n}; x_{1:n})$ be
the likelihood function ${\mathbf
P}(Y_{1:n}=q_{1:n})\prod_{i=1}^nf_{q_i}(x_i;\theta_{q_i})$ of
$q_{1:n}$, and  let $\V(x_{1:n})$ be the set of the
maximum-likelihood estimates $v(x_{1:n})\in S^n$ of $y_{1:n}$. The
elements of $\V(x_{1:n})$ are called {\em (Viterbi) alignments} and are
commonly computed by the Viterbi algorithm \cite{Viterbi, tutorial}.
If ${\mathbf P}(Y_{1:n}=q_{1:n})$ is thought of as the prior distribution 
of $Y_{1:n}$, then $v(x_{1:n})$'s also maximize the probability mass function of the posterior distribution of $Y$,
hence the term {\em maximum a posteriori (MAP) paths}. Besides their direct
significance for prediction of $Y$ from $X$, Viterbi alignments, or MAP paths, 
are also  central to the theory and applications of HMMs \cite{HMM}
in the more general setting when  any parameters of the emission
distributions $P_l$ and any of the transition probabilities $p_{ij}$, $i,j\in S$, would
also be unknown and of interest.   Therefore, asymptotic behavior of
Viterbi alignments is also crucial for the inference on the unknown
parameters \cite{HMM}, \cite{AVT4}.

To appreciate that the question of extending   $v(x_{1:n})$ {\em ad infinitum} 
is not a trivial one even if the problem of non-uniqueness of $v(x_{1:n})$ is disregarded, 
suffice it to say that an  additional  observation $x_{n+1}$
can in principle change the entire alignment based  on
$x_{1:n}$, i.e. $v(x_{1:n})$ and $v(x_{1:n+1})_{1:n}$ can disagree significantly, if not fully. 
Fortunately, the situation is not hopeless and  
in this paper we prove that in most HMMs alignments can be consistently extended 
{\em piecewise}. Specifically,  
 motifs of (contiguous) observations $z_{1:b}$, called {\it barriers}, are
observed with positive probability, 
forcing Viterbi alignments  based on extended observations $(x_{1:n},z_{1:b},x_{n+b+1:n+b+r})$, $n\ge 0$, $r\ge 1$,
to stabilize as follows:
Roughly, $v(x_{1:n}z_{1:b}x_{n+b+1:n+b+r})_{1:n}=v(x_{1:n})$ for all $x_{1:n}$ and all extensions 
$x_{n+b+1:n+b+r}$.  To be more precise, a particular state $l\in S$ and
an element $b_k$,  called a {\it node}, of the barrier $b$ can be found such
that regardless of the observations before and after $b$,
the alignment has to go through $l$ at time $u=n+k$.  
The optimality principle then insures the stabilization
$v(x_{1:n}z_{1:b}x_{n+b+1:n+b+r})_{1:u}=v(x_{1:u})$ and in particular $v_u=l$.  

Suppose now that
$x_{1:n}$ contains several barriers with nodes occurring at times 
$u_1 < \cdots < u_m\le n$. Then 
the Viterbi alignment $v(x_{1:n})$ can be
constructed piecewise as follows: Let
$v(x_{1:\infty})=(v^1,v^2,\ldots,v^m,v^{m+1})$, where $v^1$ is the
alignment based on $x_{1:u_1}$ and ending in $l$, and let 
$v^i$, for $i=2,3,\ldots,m+1$, be the conditional alignment based on 
$x_{u_{i-1}:u_i}$ given that $Y_{u_{i-1}}=l$; note that the alignments $v^i$, 
$i=2,3,\ldots,m$ also end in $l$.  Now, if a new observation
$x_{n+1}$ is added, then the last segment $v^{m+1}$ can change, but the segments 
$v^1,\ldots,v^m$ are intact. Suppose now that a realization $x_{1:\infty}$ 
contains infinitely many barriers, and hence also infinitely many nodes. Then
the {\it (piecewise) infinite} alignment $v(x_{1:\infty})$ is defined naturally 
as the infinite succession of the segments $v^1$, $v^2$, \ldots. 

In this paper, we prove that for some fixed integer $M>0$,  the probability that 
the finite random process $X_{1:M}$ generates a barrier, is positive. 
Since $X$ is  ergodic,  almost every realization  $x_{1:\infty}$
has infinitely many barriers and, therefore, the infinite piecewise 
alignment is well-defined. Apparently, the piecewise alignment
gives rise to a decoding process $v: {\cal X}^{\infty} \mapsto S^{\infty}$ via $V_{1:\infty}=v(X_{1:\infty})$, 
which we shall call  the {\it  Viterbi alignment process}. The construction
ensures that $V$ is regenerative and ergodic. Note also how
this piecewise construction naturally calls for a buffered on-line implementation
in which the memory used to store $x_{u_{i-1}:u_i}$ can be released once $v^i$ has been 
computed.    
\subsection{Previous related work and contribution of this work }\label{sec:prework}
The problem of constructing infinite Viterbi processes has been brought to 
the attention of the IEEE Information Theory community fairly recently by 
\cite{caliebe1} and \cite{caliebe2}. Although the piecewise structure of Viterbi alignments 
was already acknowledged in \cite{kogan}, to our best knowledge, the subject 
has been first seriously considered in \cite{caliebe1,caliebe2}. In these latter works, 
the existence of infinite alignments for certain special cases, such as $K=2$ and Markov 
chains with additive white Gaussian noise, has been proved. In particular, in these cases 
the authors of  \cite{caliebe1,caliebe2}  have  proved the existence of 
`meeting times' and  `meeting states', which are a special (stronger) type of nodes. 
While innovative, the main result of \cite{caliebe1} (Theorem 2) makes several restrictive assumptions and
is proved in an existential manner, which  prevents its extension beyond the $K=2$ case. 

Independently of these works, 
\cite{AVT1,AVT4,AVT3}  have  developed a  more general theory
to include the problem of estimating unknown parameters ($\theta_i$, and $p_{ij}$, $i,j\in S$). Namely, the focus of 
this theory has been the Viterbi {\em training} (VT), or {\em extraction}, algorithm \cite{jelinek0}.
Competing with EM-based procedures, this algorithm provides computationally and intuitively appealing 
estimates which, on the other hand, are biased, even in the limit when $n\to\infty$. 
In order to reduce this bias, the {\em adjusted Viterbi training} (VA) has been introduced in \cite{AVT1,AVT4,AVT3}. 
Naturally, VA relies on the existence of infinite alignments and their ergodic properties.
Although the general theory has been presented in \cite{AVT3, AVT4}, some of the main results
of the theory (Lemma 3.1 and 3.2 of \cite{AVT4}) have appeared without proof due to the limitations
of scope and size.  This paper slightly refines these results and, most importantly, presents
their complete proofs. Whereas these results are formulated for general HMMs ($K\ge 2$), 
\cite{AVTK2} has most recently considered in full detail the special case of $K=2$, generalizing
similar results of \cite{caliebe1, caliebe2}. Specifically,
it has been proved in \cite{AVTK2} that infinitely many barriers 
(and hence the infinite Viterbi alignment) 
{\it exist for any aperiodic and irreducible 2-state HMM}.  
Thus, the results presented here generalize the ones of \cite{AVTK2} and \cite{caliebe1, caliebe2}
for $K\ge 2$. It turns out that this generalization is far
from being straightforward and requires a more advanced analysis and tools. 
Furthermore, as we show below, when $K>2$, {\it not every aperiodic and
irreducible  HMM has infinitely many nodes}, undermining the piecewise
construction of infinite alignments for those models.  The disappearance of nodes 
is due to the fact that an aperiodic and irreducible Markov chain can have zeros
in the transition matrix. If this possibility is excluded, as is the case
in \cite{caliebe1, caliebe2}, 
the `meeting times' and `meeting states' of  \cite{caliebe1, caliebe2}
are sufficient to prove the existence of infinite Viterbi alignments for
many HMMs used in practice.  In their recent communication with us, the authors of
\cite{caliebe1, caliebe2} have corrected those statements in their above works
where the strict positivity of the transition matrix is implicitly assumed but
formally omitted (see \cite{AVT4} for details). 
At the same time, in order to  accommodate for zeros in the transition matrix, 
\cite{AVT4} introduced a more general notion of nodes, effectively removing 
the limitations of the notion of `meeting times' and `meeting states'.
However, the price for this generalization has been rather high due to
the interfering issue of non-uniqueness of (finite) Viterbi alignments. 
For a detailed treatment of the piecewise construction of 
the infinite alignment and process in general HMMs, and the role 
of the infinite Viterbi process for the adjusted Viterbi training theory, 
we refer to the state-of-the-art article \cite{AVT4}. 
\subsection{Organization of the rest of the paper}
In \S\ref{sec:pre} we briefly outline the construction of the
infinite alignments \S\ref{sec:process} based on \cite{AVT4}.
This includes definitions of nodes \S\ref{sec:nodes} and barriers \S\ref{sec:barriers}. Next,
\S\ref{sec:exist} states our main results which have first appeared
in \cite{AVT4} and guarantee the existence of the alignment process $V$. In
\S\ref{sec:examples}, we give a counterexample to explain the
necessity of our technical assumptions. In \S\ref{sec:proof}, we present 
a complete and detailed proof of our main results. 
This is followed in \S\ref{sec:end} by a brief discussion of
the significance of the presented results.
\section{Construction}\label{sec:pre}
\subsection{Nodes}\label{sec:nodes}
First, consider the {\it scores}
\begin{equation}\label{eq:delta}
\delta_u(l)\stackrel{\mathrm{def}}{=}\max_{q\in S^{u-1}}\Lambda\bigl((q,l); x_{1:u}\bigr).
\end{equation}
Thus, $\delta_u(l)$ is the maximum of the likelihood
of the paths terminating at $u$ in state $l$. Note that
$\delta_1(l)=\pi_lf_l(x_1)$ and the recursion below
\begin{equation*}
\delta_{u+1}(j)=\max_{l\in S}(\delta_{u}(l)p_{lj})f_j(x_{u+1})\quad\forall~u\ge 1, \forall j \in S,
\end{equation*}
helps to  verify that $\V(x_{1:n})$, the set of all the Viterbi alignments, can be written as follows: 
$\V(x_{1:n})=\left \{v\in S^n:~\forall i\in S, \delta_{n}(v_n)\ge
\delta_n(i)~\text{and}\right .$\\
$\left . \forall u:~1\le u<n,~v_u\in t(u,v_{u+1})\right \}$, where $\forall u\ge 1, \forall j \in S$,
\begin{eqnarray}
\label{viimane} 
t(u,j)\stackrel{\mathrm{def}}{=}\{l\in S: \forall i\in S~\delta_u(l)p_{lj}\ge \delta_u(i)p_{ij}\}.
\end{eqnarray}
Next, we introduce $p^{(r)}_{ij}(u)$, the maximum  of the likelihood realized along
the paths connecting states $i$ and $j$ at times $u$ and $u+r$, respectively.
Thus, $p^{(0)}_{ij}(u)\stackrel{\mathrm{def}}{=}p_{ij}$ and $\forall u\geq 1$, and $\forall r\ge 1$, let
$p^{(r)}_{ij}(u)\stackrel{\mathrm{def}}{=}\max_{q_{1:r}\in
S^r}p_{iq_1}f_{q_1}(x_{u+1})
p_{q_1q_2}f_{q_2}(x_{u+2})p_{q_2q_3}\cdots$
\begin{eqnarray}\label{eq:pr}
\cdots p_{q_{r-1}q_r}f_{q_r}(x_{u+r})p_{q_r j}.
\end{eqnarray}
Note also
\begin{align}\label{eq:prrecurse}\nonumber
\delta_{u+1}(j)&=\max_{i\in S} \bigl\{\delta_{u-r}(i) p^{(r)}_{ij}(u-r)\bigr\}f_j(x_{u+1})\quad \forall r<u,\\
p^{(r)}_{ij}(u)&=\max_{q\in S}p^{(r-1)}_{iq}(u)f_q(x_{u+r})p_{q j}.
\end{align}
\begin{definition}\label{rnode}Let  $0\le r<n$, $u\le n-r$
and let $l\in S$.  Given $x_{1:u+r}$, the first $u+r$ observations,
 $x_u$ is said to be an {\it $l$-node of order $r$} if
\begin{equation}\label{r-krit}
\delta_u(l)p^{(r)}_{lj}(u)\geq \delta_u(i)p^{(r)}_{ij}(u)\quad
\forall i,j\in S.
\end{equation}
Also, $x_u$ is said to be a node of order $r$ if it is an
 $l$-node of order $r$ for some $l\in S$; $x_u$ is said to be a
strong node of order $r$ if the inequalities in \eqref{r-krit} are strict for
every $i,j\in S, i\ne l$.\footnote{Note that if $x_u$ is a node of order $r$, it is then
also a node of any order higher than $r$. Hence, the order of a node is defined
to be the minimum such $r$.}
Let $x_{1:n}$ be such that $x_{u_i}$ is an $l_i$-node of order $r$,
$1\le i\le k$, for some $k<n$, and assume $u_{k}+r<n$ and $u_{i+1}>u_i+r$ for all $i=1,2,\ldots,k-1$.
Such  nodes are said to be {\it separated}.
\end{definition}
\subsection{Piecewise alignment}\label{sec:process}
Suppose $x_{1:n}$ is
such that  for some $u_i, r_i$, $i=1,2,\ldots, k$,
$u_1+r_1<u_2+r_2<\cdots <u_k+r_k<n$, $x_{u_i}$ is an $l_i$-node of order $r_i$.
It follows then easily from the definition of the node that there exists a
Viterbi alignment $v(x_{1:n})\in \V(x_{1:n})$ that goes through $l_i$ at
$u_i$ (i.e. $v_{u_i}=l_i$) for each $i=1,2,\ldots,k$ (see  \cite{AVT4}).
It is not difficult to verify that such $v(x_{1:n})$ can actually be
computed as follows: Obtain $v^1$, a path that is optimal
among all those that end at $u_1$ in $l_1$. (Note that unless the order of
the node $x_{u_1}$ is 0, $v^1$ need not be in $\V(x_{1:u_1})$.)
Given $x_{u_1+1:u_2}$, continue on by taking $v^2$ to
be a maximum likelihood path from $l_1$ to $l_2$.  That is, 
$v^2$ maximizes the constrained likelihood under the initial distribution 
$(p_{l_1\cdot})$ and the constraint $v^2_{u_2-u_1}=l_2$. 
Now, $(v^1,v^2)$ maximizes the likelihood given $x_{1:u_2}$ over all paths ending with $l_2$.
Similarly, we define the pieces $v^3,\ldots,v^k$. Finally, $v^{k+1}$ is chosen to
maximize the (unconstrained) likelihood given $x_{u_{k+1}:n}$
under the initial distribution $(p_{l_k\cdot})$.

The separated nodes assumption $u_{i+1}>u_i+r$,
$1\le i<k$,  is not restrictive at all  since
it is always possible to choose  from any infinite sequence of nodes
an infinite subsequence of separated ones. The reason for this requirement has to do
with the non-uniqueness of alignments and is as follows. 
The fact that $x_{u_i}$ is an $r$th order $l_i$-node  
guarantees that when backtracking from $u_{i}+r$ down to $u_i$, ties
(if any)
can be broken in such a way that, regardless of the values of  $x_{u_i+r+1:n}$ and how
ties are broken in between $n$ and $u_i+r$, the
alignment goes through $l_i$ at $u_i$. At the same time, segment $u_i,\ldots,u_i+r$ is `delicate',
that is, unless $x_{u_i}$ is a strong node,
breaking the ties arbitrarily within $u_i,\ldots,u_i+r$ can result in  $v_{u_i}\neq l_i$.
Hence, when neither $x_{u_i}$ nor
$x_{u_{i+1}}$ is strong and 
$u_{i+1}\leq u_{i}+r$, breaking the ties in favor of
$x_{u_i}$ can result in $v_{u_{i+1}}\neq l_{i+1}$. Clearly,  such a pathological situation is impossible if $r=0$
and might also be rare in practice even for $r>0$.


To formalize the piecewise construction, let
\begin{align*}
&\W^{l}(x_{1:n})\stackrel{\mathrm{def}}{=}
\{v\in S^n:~v_n=l\\
& \Lambda(v; x_{1:n})\ge \Lambda(w; x_{1:n})\quad \forall w\in S^n: w_n=l\},
\end{align*}
$\V^{l}(x_{1:n})\stackrel{\mathrm{def}}{=}\{v\in \V(x_{1:n}): v_n=l\}$ be the set  of  
maximizers of the constrained  likelihood,
and the subset of  maximizers of the (unconstrained) likelihood,
respectively, all elements of which  go through $l$ at $n$. Note
that unlike $\W^{l}(x_{1:n})$, $\V^{l}(x_{1:n})$ might be empty.
It can be shown that $\V^{l}(x_{1:n})\neq \emptyset\Rightarrow
\V^{l}(x_{1:n})=\W^{l}(x_{1:n})$.
Also, let  subscript the $(l)$ in $\W_{(l)}^{m}(x_{1:n})$ and $\V_{(l)}(x_{1:n})$
refer to  $(p_{li})_{i\in S}$ being used as  the initial distribution in place of $\pi$.
With these notations, the piecewise alignment is
$v=(v^1,\ldots,v^{k+1})\in \V(x_{1:n})$, where
\begin{align}\label{r-piecesI}
v^1\in &\W^{l_1}(x_{1:u_1}),\quad v^{k+1}\in\V_{(l_{k})}(x_{u_{k}+1:u_n})\nonumber \\
v^i\in &\W_{(l_{i-1})}^{l_i}(x_{u_{i-1}+1:u_i}),~2\le i\le k.
\end{align}
Moreover, for $i=1,2,\ldots,k$, the partial paths
$w(i)\stackrel{\mathrm{def}}{=}(v^1,\ldots,v^i)\in \W^{l_i}(x_{1:u_i})$.

If $x_{1:\infty}$ has infinitely many (separated) nodes $\{x_{u_k}\}_{k\ge 1}$
then $v( x_{1:\infty})$, an {\it infinite piecewise alignment based 
on the node times $\{u_k(x_{1:\infty})\}_{k\ge 1}$} can be defined as follows:  
If  the sets $\W_{(l_{i-1})}^{l_i}(x_{u_{i-1}+1:u_i})$,
$i=2,\ldots, k$ as well as $\V_{(l_k)}(x_{u_k+1:n})$ and
 $\W^{l_1}(u_1,x_{1:u_1})$ are singletons, then  \eqref{r-piecesI} immediately defines a
unique  infinite alignment  $v(x_{1:\infty})=(v^1(x_{1:u_1}),
v^2(x_{u_1+1:u_2}),\ldots)$. Otherwise, ties must be broken.
If we want our infinite alignment process $V$ to be regenerative (see \cite{AVT4}), 
a  natural consistency condition must be imposed on
rules to  select unique $v(x_{1:n})$  from $\W^{l_1}(x_{1:u_1})\times 
\W_{(l_1)}^{l_2}(x_{u_1+1:u_2})\times\cdots\times\W_{(l_{k-1})}^{l_k}(x_{u_{k-1}+1:u_k})\times 
\V_{(l_{k})}(x_{u_{k}+1:n})$. In \cite{AVT4}, resulting infinite
alignments, as well as decoding $v:~\mathcal{X}^\infty \to S^\infty$ based on such 
alignments,  are  called {\em proper}. This condition is, perhaps, best understood
by the following example. Suppose for some $x_{1:5}\in \mathcal{X}^5$, 
$\W_{(1)}^1(x_{1:5})=\{12211, 11211\}$, and suppose the tie is broken in 
favor of $11211$. Now, whenever 
$\W_{(l)}^1(x'_{1:4})$ contains $\{1221, 1121\}$, we naturally require that 
$1221$   not be selected. In particular, 
we select $1121$ from $\W_{(1)}^1(x_{1:4})=\{1221,
1121\}$. Subsequently, 
$112$ is selected from $\W_{(1)}^2(x_{1:3})=\{122, 112\}$, and so on.
{\em It can be shown that a decoding by piecewise alignment \eqref{r-piecesI} with ties 
broken in favor of min (or max) under the reverse lexicographic ordering of $S^n$, 
$n\in \mathbb{N}$, is a proper decoding.} 

Note also that we break ties locally, i.e. within 
individual intervals $u_{i-1}+1,\ldots,u_i$, $i\ge 2$, enclosed by 
adjacent nodes. This is in contrast to global ordering of $\V(x_{1:n})$,
such as the one in \cite{caliebe1, caliebe2}. Since a global order 
need not respect  decomposition \eqref{r-piecesI},
it can fail to produce an infinite alignment  going through infinitely
many nodes unless the nodes are strong.

\subsection{Barriers}\label{sec:barriers}
Recall (Definition \ref{rnode}) that nodes of order $r$ at time $u$ are defined {\em  relative} to the entire
realization $x_{1:u+r}$. Thus, whether $x_u$ is a node or not depends,
in principle, on all observations up to $x_u$.

We show below that typically a block $x^b_{1:k}\in \mathcal{X}^k$ ($k\ge r$) can be found such 
that for any $w\ge 1$ and for any $x'_{1:w}\in \mathcal{X}^w$, $(w+k-r)$th element of 
$(x'_{1:w},x^b_{1:k})$ is a node of order $r$ (relative to $(x'_{1:w},x^b_{1:k})$).  
Sequences $x^b_{1:k}$ that ensure existence of such persistent nodes are called {\it barriers} in
\cite{AVT4}. Specifically,
\begin{definition}\label{def:barrier}
Given $l\in S$, $x^b_{1:k}\in {\cal X}^k$ is called an (strong)
$l$-{\it barrier} of order $r\ge 0$ and length $k\ge 1$ 
if,  for any $w\ge 1$ and for every $x'_{1:w}\in \mathcal{X}^w$, 
$(x'_{1:w},x^b_{1:k})$ is such that  $(x'_{1:w},x^b_{1:k})_{w+k-r}$
is an (strong) $l$-node of order $r$. 
\end{definition}
\section{Existence}\label{sec:exist}
\subsection{Clusters and main results}\label{sec:clusters}
For each $i\in S$, let $$G_i\stackrel{\mathrm{def}}{=}\{x\in\mathcal{X}: f_i(x)>0\}.$$
\begin{definition}We call a subset $C\subset S$  {\it a cluster} if the following conditions are satisfied:
$$\min_{j\in C}P_j(\cap _{i\in C}G_i)>0,~{\rm and}\,\max_{j\not\in C}P_j(\cap _{i\in C}G_i)=0.$$
\end{definition}
Hence, a cluster is a maximal subset of states such that $G_C=\cap _{i\in C}G_i$, the intersection
of the supports of the corresponding emission distributions, is  `detectable'.
Distinct clusters need not be disjoint and  a
cluster can consist of a single state. In this latter case such a state is not
hidden, since it is  exposed by any observation it emits.
When $K=2$, $S$ is the only  cluster possible,
since otherwise all observations would expose their states and
the underlying Markov chain would cease to be hidden.
In practice,  many other HMMs  have  the entirety of $S$ as their (necessarily unique) cluster.

We now state the main results. For every state $l\in S$, let
\begin{equation}\label{eq:pmaxtrans}
p^*_{l}=\max_{j}p_{jl}.
\end{equation}
\begin{lemma}\label{neljas} Assume that for each state $l\in S$,
\begin{equation}\label{lll}
P_l\left(\left\{x\in\mathcal{X}:~f_l(x)p^*_{l}> \max_{i,i\ne
l}f_i(x)p^*_{i}\right\}\right)>0.
\end{equation}
Moreover, assume that there exists a cluster $C\subset S$ and a positive
integer $m$ such that the $m$th power of the sub-stochastic
matrix $\mathbb{Q}=(p_{ij})_{i,j\in C}$ is strictly positive.
Then, for some integers $M$ and $r$, $M>r\geq 0$,
there exist a set $B=B_1\times \cdots \times B_M \subset {\cal X}^M$, an $M$-tuple
of states $q_{1:M}\in S^M$ and a state $l\in S$, such that
every  $x_{1:M}\in B$ is an
$l$-barrier of order $r$ (and length $M$), $q_{M-r}=l$ and
\begin{eqnarray}\label{eq:result}\nonumber
\P\left(X_{1:M}\in B, \quad Y_{1:M}=q_{1:M}\right)>0.
\end{eqnarray}
\end{lemma}

Lemma \ref{neljas} implies that $\P(X_{1:M}\in B)>0$. Also, since every element of 
$B$ is a
barrier of order $r$, the ergodicity of $X$ therefore guarantees that {\em almost every
realization of $X$
contains infinitely many  $l$-barriers of order $r$.} {\em Hence, almost every realization of $X$
also has infinitely many $l$-nodes of order $r$.}

In two state HMMs, $S$ is the only cluster (otherwise the Markov chain would not be hidden), hence
$\mathbb{Q}=\mathbb{P}$. The irreducibility and aperiodicity in this case imply strict positivity of $\mathbb{P}^2$.
Thus, the only condition to be verified is \eqref{lll}, which in this case writes as 
$P_1\left(\left\{x\in \mathcal{X}:~f_1(x)p^*_{1}> f_2(x)p^*_{2}\right\}\right)>0$
and  $P_2\left(\left\{x\in \mathcal{X}:~f_2(x)p^*_{2}> f_1(x)p^*_{1}\right\}\right)>0$. 
In \cite{AVTK2}, it is shown that in the case of two state HMMs,  one of these two positivity conditions is always met, 
which, in fact, turns out to be sufficient for the existence of infinitely many strong barriers in this ($K=2$) case. 
Thus, {\em any two state HMM with
irreducible and aperiodic $Y$ has infinitely many strong barriers.}  
Lemma  \ref{neljas} significantly generalizes this and associated results of \cite{AVTK2}. The case
$K=2$ is special in several respects, hence the generalization is technically involved, and in particular the
CLT-based proof of the existence of infinitely many nodes in \cite{caliebe1} (Theorem 2) does not apply 
when $K>2$.

For certain technical reasons, instead of extracting subsequences of separated
nodes from general infinite sequences of nodes guaranteed by Lemma \ref{neljas},
we achieve node separation by adjusting the notion of barriers.
Namely, note that two $r$th-order $l$-barriers $x_{j:j+M-1}$ and $x_{i:i+M-1}$
might be in $B$ with $j<i\leq j+r$, implying that the associated nodes
$x_{j+M-r-1}$ and $x_{i+M-r-1}$ are not separated. 
Thus, we impose on  $B$ 
the following  condition:
\begin{equation}\label{def:separated}
x_{j:{j+M-1}},x_{i:{i+M-1}}\in B,\, i\ne j~
\Rightarrow |i-j|>r.\end{equation} If \eqref{def:separated} holds,
we say that the barriers from $B\subset {\cal X}^M$ are {\it
separated}. This is often easy to achieve by a simple extension of $B$
as shown in the following example. Suppose
there exists $x\in {\cal X}$ such that $x\not \in B_m$, for all $m=1,2,\ldots,M$.
All elements of $B^*\stackrel{\mathrm{def}}{=}\{x\}\times B$ are
evidently barriers, and moreover, they are now  separated. The
following Lemma incorporates a more general version of the above example.
\begin{lemma}\label{separated}
Suppose  the assumptions of Lemma \ref{neljas} are satisfied. Then,
for some integers $M$ and $r$, $M>r\geq 0$, there exist 
$B=B_1\times \cdots \times B_M \subset {\cal X}^M$, $q_{1:M}\in S^M$, and  $l\in S$,
such that every  $x^b_{1:M}\in B$ is a separated  $l$-barrier
of order $r$ (and length $M$), $q_{M-r}=l$, and
${\bf P}\left(X_{1:M}\in B, \quad Y_{1:M}=q_{1:M}\right)>0$.
\end{lemma}
\subsection{Counterexamples}\label{sec:examples}
The condition on $C$ in Lemma \ref{neljas} might seem technical
and even unnecessary. We next give an example of an HMM where
the cluster condition is not met and no node (barrier) can occur.
Then, we will modify the example
to enforce the cluster condition and consequently gain barriers.
\begin{example}\label{ex:2.3} Let $K=4$ and consider an ergodic Markov chain
with transition matrix
$$\mathbb{P}=
\left(%
\begin{array}{cccc}
  {1\over 2} & 0 & 0 &  {1\over 2} \\
  0 &  {1\over 2} &  {1\over 2} & 0 \\
   {1\over 2} & 0 &  {1\over 2} & 0 \\
  0 &  {1\over 2} & 0 &  {1\over 2}\\
\end{array}%
\right).
$$
Let the emission distributions be such that \eqref{lll} is
satisfied and $G_1=G_2$ and $G_3=G_4$ and $G_1\cap G_3=\emptyset$.
Hence, in this case there are two disjoint clusters $C_1=\{1,2\}$,
$C_2=\{3,4\}$. The matrices $\mathbb{Q}_i$  corresponding to $C_i$,
$i=1,2$ are
$$\mathbb{Q}_1=\mathbb{Q}_2=
\left(%
\begin{array}{cc}
   {1\over 2} & 0 \\
  0 &  {1\over 2} \\
\end{array}%
\right).$$ Evidently, the cluster assumption of Lemma \ref{neljas} is not
satisfied. Note also that the alignment cannot change (in one step) its state to
the opposite one within the same cluster. Since the  supports $G_{1,2}$ and $G_{3,4}$ are disjoint, 
any observation  exposes the corresponding cluster. Hence any sequence of observations can
be regarded as a sequence of  blocks emitted from
alternating clusters. However, the alignment inside each block stays constant.
It can be shown that in this case no $x_u$ can be a node (of any order) 
for any $n>1$, $x_{1:n}\in\mathcal{X}^n$, and $1\le u<n$. 
\end{example}

Let us modify the HMM in Example \ref{ex:2.3} to ensure
the assumptions of Lemma \ref{neljas}. 
\begin{example}Let $\epsilon$ be such that $0<\epsilon<{1\over 2}$ and let us replace $\mathbb{P}$ by the following
transition matrix $$
\left(%
\begin{array}{cccc}
  {1\over 2}-\epsilon & \epsilon & 0 &  {1\over 2} \\
  \epsilon &  {1\over 2}-\epsilon &  {1\over 2} & 0 \\
   {1\over 2} & 0 &  {1\over 2} & 0 \\
  0 &  {1\over 2} & 0 &  {1\over 2}\\
\end{array}%
\right).
$$
Let the emission distributions be as in the previous example. In
this case, the cluster $C_1$ satisfies the assumption of Lemma
\ref{neljas}. As previously, every observation exposes its
cluster. Lemma \ref{neljas} now applies to guarantee barriers and nodes. 
To be more specific, let $\epsilon=1/4$, $f_1(x)=\exp(-x)_{x\ge 0}$, 
$f_2(x)=2\exp(-2x)_{x\ge 0}$, and $f_3(x)=\exp(x)_{x\le 0}$, 
$f_4(x)=2\exp(2x)_{x\le 0}$. It can then be verified that
if $x_{1:2}=(1,1)$ then $x_1$ is a $1$-node of order 2. Indeed, in
that case any element of  $B=(0,+\infty)\times (\log(2), +\infty)\times (0,+\infty)$ is
a $1$-barrier of order 2. 
\end{example}

Another way to modify the HMM in  Example \ref{ex:2.3}
to enforce the assumptions of Lemma \ref{neljas} is to change
the emission probabilities. Namely, assume that the supports $G_i$,
$i=1,\ldots,4$ are such that $P_j(\cap_{i=1}^4 G_i)>0$ for all $j\in S$, and
\eqref{lll} holds. Now, $S=\{1,\ldots, 4\}$ is the only cluster.
Since the matrix $\mathbb{P}^2$ has all its entries positive, the conditions
of Lemma \ref{neljas} are now satisfied and barriers can now be constructed. 

\section{Proof of the main result}\label{sec:proof}
\subsection{Proof of Lemma \ref{neljas}}\label{sec:proofneljas}
The proof below is a rather direct construction which is, however, technically
involved. In order to facilitate the exposition of this proof, we have divided it into 17
short parts as follows.
\subsubsection{$\mathcal{X}_l\subset \mathcal{X}$}
\label{subsec:Xl}
It follows from the assumption \eqref{lll} and finiteness of $S$ that there exists an
$\epsilon>0$ such that for all $l\in S$ $P_l({\cal X}_l)>0,$ where
\begin{eqnarray}\label{tarn}
{\cal X}_l\stackrel{\mathrm{def}}{=}\Big\{x\in\mathcal{X}: \max_{i,i\ne
l}p^*_i f_i(x)<(1-\epsilon)p^*_lf_l(x)\Bigr\}.&&
\end{eqnarray}
(Note that $p^*_l>0$ for all $l\in S$ by irreducibility of $Y$.) Also note that
 ${\cal X}_l,l\in S$ are disjoint and have positive reference measure
$\lambda({\cal X}_l)>0$.
%
\subsubsection{$\mathcal{Z}\subset\mathcal{X}$ and $\delta-K$ bounds on cluster densities $f_i$, $i\in C$}
\label{subsec:Z}
Let $C$ be a cluster as in the assumptions of the Lemma. The existence of $C$ implies the
existence of a set $\hat{{\cal Z}}\subset \cap_{i\in C }G_i$ and $\delta>0$,
 such that $\lambda (\hat{{\cal Z}})>0$, and $\forall z\in
{\hat{\cal Z}}$, the following statements hold:
\begin{enumerate}[(i)]
    \item $\min_{i\in C} f_i(z)>\delta$;
    \item $\max_{j\not \in C} f_j(z)=0$.
\end{enumerate}
Indeed, $\min_{j\in C}P_j(\cap _{i\in C}G_i)>0$ implies (and indeed is equivalent to) $\lambda(\cap _{i\in C}G_i)>0$. 
The latter implies the existence of $\hat{\cal Z}\subset \cap _{i\in C}G_i$ with positive $\lambda$-measure 
and  $\delta>0$ such that (i) holds. Since
$\lambda(\cap _{i\in C}G_i)>0$, the condition $P_j(\cap _{i\in C}G_i)=0$ for $j\not \in C$  implies (is equivalent to)
$f_j=0$  $\lambda$-{\em almost everywhere} on $\cap _{i\in C}G_i$. Thus, 
$\max_{j\not \in C} f_j=0$  $\lambda$-{\em almost everywhere} on $\cap _{i\in C}G_i$,  which implies (ii).

Evidently, $K>0$ can be chosen sufficiently large to make
$\lambda (\{z\in\mathcal{X}:~f_i(z)\ge K\})$ arbitrarily small, and in particular, to
guarantee that
$\lambda (\{z\in\mathcal{X}:~f_i(z)\ge K\})<\frac{\lambda (\hat{{\cal Z}})}{|C|}$,
where $|C|$ is the size of $C$.
Clearly then, redefining
$\hat{{\cal Z}}\stackrel{\mathrm{def}}{=}
\hat{{\cal Z}}\cap\{z\in\mathcal{X}:~f_i(z)<K,~i\in C\}$ preserves
$\lambda (\hat{{\cal Z}})>0$.
Next, consider
\begin{equation}\label{moot}
\lambda(\hat{{\cal Z}}\backslash (\cup_{l\in S}{\cal X}_l)).
\end{equation}
If \eqref{moot} is positive, then define
\begin{equation}
  \label{eq:defineZ1}
{\cal Z}\stackrel{\mathrm{def}}{=}\hat{{\cal Z}}\backslash (\cup_{l\in S}{\cal X}_l).
\end{equation}
If \eqref{moot} is zero, then there must be  $s\in C$ such that
$$\lambda(\hat{{\cal Z}}\cap {\cal X}_s)>0$$
and in this case, let
\begin{equation}
  \label{eq:defineZ2}
{\cal Z}\stackrel{\mathrm{def}}{=}\hat{{\cal Z}}\cap {\cal X}_s.
\end{equation}
Such $s\in S$ must clearly exist since $\lambda(\hat{{\cal Z}})>0$ but
$\lambda(\hat{{\cal Z}}\backslash (\cup_{l\in S}{\cal X}_l))=0$.
To see that $s$ must necessarily be in the cluster $C$,
note  
 $\forall s\not\in C$, $f_s(z)=0$ $\forall z\in{\hat{\cal Z}}$, which implies
${\hat{\cal Z}}\cap {\cal X}_s =\emptyset$.
\subsubsection{Sequences $\mathbf{s}$,  $\mathbf{a}$, and $\mathbf{b}$ of  states in $S$}
\label{subsec:abs} Let us define an auxiliary sequence of states $q_1$, $q_2$, and so on, as follows:
If \eqref{moot} is zero, that is, if
${\cal Z}=\hat{{\cal Z}}\cap {\cal X}_s$ for some $s\in C$, then define $q_1=s$,
otherwise let $q_1$ be an arbitrary state in $C$.
Let $q_2$ be a state with maximal probability of transition to $q_1$, i.e.:
$p_{q_2\,q_1}=p^*_{q_1}$
Suppose $q_2\ne q_1$. Then find $q_3$ with $p_{q_3\,q_2}=p^*_{q_2}$.
If $q_3\not\in \{q_1,q_2\}$, find $q_4:~p_{q_4\,q_3}=p^*_{q_3}$, and so on.
Let $U$ be the first index such that $q_U\in \{q_1,\ldots,q_{U-1}\}$,
that is, $q_U=q_T$ for some $T<U$.
This means that there exists a sequence of states  $\{q_T,\ldots, q_U\}$ such that
 \begin{itemize}
    \item $q_T=q_{U}$
    \item $q_{T+i}=\arg\max_{j} p_{jq_{T+i-1}},\quad i=1, \ldots, U-T.$
\end{itemize}
To simplify the notation and without loss of generality, assume $q_U=1$.
Reorder and rename the states as follows:
\begin{align*}
&s_1\stackrel{\mathrm{def}}{=}q_{U-1},\, s_2\stackrel{\mathrm{def}}{=}q_{U-2},\ldots,
s_i\stackrel{\mathrm{def}}{=}q_{U-i},\ldots,\nonumber \\
&s_L\stackrel{\mathrm{def}}{=}q_T=1\quad i=1,\ldots,L\stackrel{\mathrm{def}}{=}U-T,\\
&
a_1\stackrel{\mathrm{def}}{=}q_{T-1},\,a_2\stackrel{\mathrm{def}}{=}q_{T-2},\ldots,a_P\stackrel{\mathrm{def}}{=}q_1,
\end{align*}
where $P\stackrel{\mathrm{def}}{=}T-1$. Hence,
\begin{align*}
&\{q_1,\ldots,q_{T-1},q_T,q_{T+1},\ldots,q_{U-1},q_U\}=\\
&\{a_P,\ldots,a_1,1,s_{L-1},\ldots,s_1,1\}.
\end{align*}
Note that if $T=1$, then $P=0$ and
$\{q_1,\ldots,\ldots,q_{U-1},q_U\}=\{1,s_{L-1},\ldots,s_1,1\}.$
We have thus introduced special sequences $\mathbf{a}=(a_1,a_2,\ldots,a_P)$ and
$\mathbf{s}=(s_1,s_2,\ldots,s_{L-1},1)$.
Clearly,
\begin{align}\label{klass}
p_{s_{i-1}\,s_i}=&p^*_{s_i},~ i=2,\ldots, L,~ p^*_{s_1}=p_{1\,s_1}\nonumber \\
p_{a_{i-1}\,a_i}=&p^*_{a_i},~ i=2,\ldots, P,~ p^*_{a_1}=s_L=1.
\end{align}
Next, we are going to exhibit $\mathbf{b}=(b_1,\ldots,b_R)$, another auxiliary sequence for some $R\ge 1$,
characterized as follows:
\begin{enumerate}[(i)]\label{eq:defineb}
    \item $b_R=1$;
    \item $\exists$ $b_0\in C$ such that
    $p_{b_0\,b_1}p_{b_1\,b_2}\cdots
    p_{b_{R-1}\,b_R}>0$;
    \item if $R>1$, then $b_{i-1}\ne b_i$ for every $i=1,\ldots,R$.
\end{enumerate}
Thus, the path $b_{1:R}$ connects cluster $C$ to state 1
in $R$ steps. Let us also require that $R$ be minimum such. Clearly
such $\mathbf{b}$ and $b_0$ do exist due to irreducibility of $Y$. Note also that minimality
of $R$ guarantees (iii) (in the special case of $R=1$ it may happen that $b_1=1\in S$
and $p_{1\,1}>0$, in which case $b_0$ can be taken to be also $1$).
\subsubsection{Determining $k$}
\label{subsec:k}
Let $\mathbb{Q}^{m}$ be the $m$th power of  the
sub-stochastic matrix $\mathbb{Q}=(p_{ij})_{i,j\in C}$; let $q_{ij}$ be the entries of $\mathbb{Q}^{m}$.
By the hypothesis of the Lemma, $q_{ij}>0$ $\forall i,j\in C$. This means that for every $i,j\in C$,
there exists a positive probability
path from $i$ to $j$ of length $m$.
Let $q^*_{ij}$ be the probability of a maximum probability path
from $i$ to $j$. In other words, for every $i,j\in C$, there exist
states $w_1,\ldots,w_{m-1}\in C$ such that
\begin{eqnarray}
\label{maxprobpath}
p_{iw_1}p_{w_1w_2}\cdots
p_{w_{m-1}w_{m-1}}p_{w_{m-1}j}=q^*_{ij}>0.&&%
\end{eqnarray}Let us define
\begin{eqnarray}
\label{eq:qpositive}
 q&=&\min_{i,j\in C}q^*_{i\,j}>0,\quad\text{and}\\
\label{eq:maxtransandratio}  
 A&=&\max_{i\in S}\max_{j\in S}\left\{{p^*_i\over p_{ji}}:p_{ji}>0\right\},
\end{eqnarray}
 where $p^*_i$'s are as defined in \eqref{eq:pmaxtrans}.
Choose $k$ sufficiently large for the following to hold:
\begin{equation}
  \label{eq:largek}
 (1-\epsilon)^{k-1}<q^2\left({\delta\over K}\right)^{2m}A^{-R},
\end{equation}
where $\epsilon$  is as in \eqref{tarn}
and $\delta$ and $K$ are as introduced in \S\ref{subsec:Z}.
\subsubsection{The $s$-path}
\label{subsec:spath}
We now fix the state sequence
\begin{equation}\label{path}
b_0,b_1,\ldots, b_R,s_1,s_2,\ldots, s_{2Lk},a_1,\ldots,a_P,
\end{equation}
where $s_{Lj+i}=s_i$, $j=1,\ldots,2k-1$, $i=1,\ldots,L$,
(and in particular $s_{Lj}=1$, $j=1,\ldots,2k$). The sequence \eqref{path} will be called
the {\em $s$-path}. The $s$-path is a concatenation  of $2k$ $\mathbf{s}$ cycles
$s_{1:L}$, the beginning  and the end of which are connected to the cluster $C$
via positive probability paths $\mathbf{b}$ and $\mathbf{a}$, respectively
(recall that $a_P=q_1\in C$ and $b_R=1$ by construction). Additionally, the
$b_R,s_1,s_2,\ldots, s_{2Lk},a_1,\ldots,a_P$-segment of the $s$-path
\eqref{path} has the important property \eqref{klass}, i.e. every consecutive transition along
this segment occurs with the maximal transition probability given
its destination state. 
(However, $\mathbf{b}$, the beginning of the $s$-path, need  not satisfy
this property.) The $s$-path is almost ready to serve as $q_{1:M}$ promised by the Lemma
and its conversion to $q_{1:M}$ will be completed in \S\ref{subsec:M}.
In fact, the  idea of the Lemma and its proof is to exhibit (a cylinder subset of) observations such that once
emitted along the $s$-path, these observations would trap the Viterbi backtracking so that the latter winds up
on the $s$-path.
That will guarantee that an observation corresponding to the beginning of the $s$-path, is a node.
\subsubsection{The barrier}
\label{subsec:barrier}
Consider the following sequence of observations
\begin{align}\nonumber
&z_0,z_1,\ldots,z_{m},y'_1,\ldots,y'_{R-1},y_0,y_1,\ldots,y_{2Lk},\\
\label{blokk}
&y^{''}_1,\ldots,y^{''}_P,z'_1,\ldots,z'_{m},\end{align}
where 
\begin{align}\nonumber
&z_0,z_i,z'_i\in {\cal Z},\quad i=1,\ldots,m;\\\nonumber
&y'_i\in  {\cal X}_{b_i},\quad i=1,\ldots,R-1;\\\nonumber
&y_0\in {\cal X}_1,\quad y_{i+Lj}\in {\cal X}_{s_i},~ j=1,\ldots,2k-1,i=1,\ldots,L\\\nonumber
&y^{''}_i\in {\cal X}_{a_i},\quad i=1,\ldots,P.
\end{align}
From this point on throughout \S\ref{subsec:onenode}, we shall be proving that $y_{Lk}$ is a
1-node of order $(kL+m+P)$, and, therefore, that \eqref{blokk} is a 1-barrier of order $(kL+m+P)$.

First, let $u\ge 2Lk+2m+1+P+R$ and let $x_{1:u}$ be any sequence of observations 
containing
the sequence \eqref{blokk}
in the tail.
\subsubsection{$\alpha$, $\beta$, $\gamma$, $\eta$}
\label{subsec:alpha-eta}
Recall the definition of the scores $\delta_u(i)$ \eqref{eq:delta} and the maximum partial
likelihoods $p^{(r)}_{i\,j}(u)$ \eqref{eq:pr}.
Now, we need to introduce the following abbreviated notation. 
For any $i,j\in S$ and appropriate $r\ge 0$, let 
\begin{align}\label{eq:shortscores}
\delta_i(y_l)&\stackrel{\mathrm{def}}{=}\delta_{u-P-m-2kL+l}(i)\quad \forall l: 0\le l\le 2kL\nonumber\\
p_{ij}^{(r)}(y_l)&\stackrel{\mathrm{def}}{=}p_{ij}^{(r)}(u-P-m-2kL+l),\\
\nonumber
p_{ij}^{(r)}(y'_l)&\stackrel{\mathrm{def}}{=}p_{ij}^{(r)}(u-P-m-2kL-R+l) \quad \forall l:\\
\nonumber
 & 1\le l\le R-1,\\
\nonumber
 \delta_i(z_l)&\stackrel{\mathrm{def}}{=}\delta_{u-2Lk-2m-P-R+l}(i)\quad \forall l:0\le l\le m,\\ \nonumber
p_{ij}^{(r)}(z_l)&\stackrel{\mathrm{def}}{=}p_{ij}^{(r)}(u- 2Lk-2m-P-R+l),
\end{align}
\begin{align}
\nonumber
\delta_i(z'_l)&\stackrel{\mathrm{def}}{=}\delta_{u-m+l}(i)\quad \forall l: 1\le l\le m,\nonumber \\
p_{ij}^{(r)}(z'_l)&\stackrel{\mathrm{def}}{=}p_{ij}^{(r)}(u-m+l).
\end{align}
Also, we will be frequently using the scores corresponding to $z_0$,
$y'_1$, $y_{Lk}$, and $y_{2Lk}$, hence the following further abbreviations:
$$\alpha_i\stackrel{\mathrm{def}}{=}\delta_i(z_0),~ \beta_i\stackrel{\mathrm{def}}{=}\delta_i(z_m),~
\gamma_i\stackrel{\mathrm{def}}{=}\delta_i(y_0),~\eta_i\stackrel{\mathrm{def}}{=}\delta_i(y_{Lk}).$$
Note  that $\forall j\not \in C$, $f(z_0)=f_j(z'_l)=f_j(z_l)=0$, $l=1,\ldots,m$
by construction of ${\cal Z}$ (\S\ref{subsec:Z}). Hence,
$\alpha_j=\beta_j=0$ $\forall j\not\in C$, and a more general implication is
that for every $j\in S$
\begin{align}\label{kreeka1}
\beta_j&=\max_{i\in
C}\alpha_ip_{ij}^{(m-1)}(z_0)f_j(z_m) \\
       &=\alpha_{i_{\beta}(j)}p_{i_{\beta}(j)\,j}^{(m-1)}(z_0)f_j(z_m)~
\text{for some }i_{\beta}(j)\in C;\nonumber \\
\label{kreeka2} \gamma_j&=\max_{i\in C}
\beta_ip_{ij}^{(R-1)}(z_m)f_j(y_0) \\
       &=\beta_{i_{\gamma}(j)}p_{i_{\gamma}(j)\,j}^{(R-1)}(z_m)f_j(y_0)~
\text{for some }i_{\gamma}(j) \in C.\nonumber
\end{align}
Also, we will use the following representation of $\eta_j$ in terms of $\gamma$:
\begin{align}
\label{kreeka3} \eta_j&=\max_{i\in S}
\gamma_ip_{i\,j}^{(kL-1)}(y_0)f_j(y_{kL})\\
            &=\gamma_{i_{\eta}(j)}p_{{i_{\eta}(j)}\,j}^{(kL-1)}(y_0)f_j(y_{kL})~
\text{for some }i_{\eta}(j)\in S.\nonumber
\end{align}
\subsubsection{Bounds on $\beta$}
\label{subsec:betabound}
Recall (\S\ref{subsec:abs}) that $b_0\in C$. We show that for every $j\in S$
\begin{align}\label{betas}
\beta_j <q^{-1}\Bigl({K\over \delta}\Bigr)^{m}\beta_{b_0}.
\end{align}
Fix $j\in S$ and consider $\alpha_{i_\beta(j)}$ from \eqref{kreeka1}. Let
$v_1,\ldots,v_{m-1}$ be a path that realizes $p_{ij}^{(m-1)}(z_0)$.
\\Then
$\beta_j=\alpha_{i_\beta(j)}p_{i_\beta(j)\,v_1}f_{v_1}(z_1)p_{v_1\,v_2}f_{v_2}(z_2)\cdots$
$p_{v_{m-1}\,j}f_{j}(z_m)<\alpha_{i_\beta(j)}K^m.$ (The last inequality follows
from \eqref{eq:defineZ1}, \eqref{eq:defineZ2}.)  Let $w_1,\ldots, w_{m-1}$ be a
maximum probability path from $i_{\beta(j)}$ to $b_0$ as in \eqref{maxprobpath}. Thus,
\begin{align*}
&\beta_{b_0}&\geq& \alpha_{i_\beta(j)}p_{i_\beta(j)\,b_0}^{(m-1)}(z_0)f_{b_0}(z_m)\\
&&\geq&
\alpha_{i_\beta(j)}p_{i_\beta(j)\,w_1}f_{w_1}(z_1)p_{w_1\,w_2}f_{w_2}(z_2)\cdots\\
&&&\cdots p_{w_{m-1}\,b_0}f_{b_0}(z_m)
\geq \alpha_{i_\beta(j)} q \delta^{m}.
\end{align*}
 (The last inequality again
follows from \eqref{eq:defineZ1}, \eqref{eq:defineZ2}.) 
Since $q>0$ \eqref{eq:qpositive}, we thus obtain:
\begin{equation*}
\beta_j<\alpha_{i_\beta(j)}K^m\leq \frac{\beta_{b_0}}{q \delta^{m}} K^m,
\end{equation*}
as required.
\subsubsection{Likelihood ratio bounds}
\label{subsec:llbound}
We next prove the following claims
\begin{align}\label{krim}
    p^{(L-1)}_{i1}(y_{lL})&\leq p^{(L-1)}_{11}(y_{lL})\nonumber \\
    \forall i\in S&\quad \forall l=0,\ldots,2k-1,\\
    \label{tsiv}
    {p^{(L-1)}_{ij}(y_{lL})f_j(y_{(l+1)L})\over
    p^{(L-1)}_{11}(y_{lL})f_1(y_{(l+1)L})}&<1-\epsilon\nonumber \\
   \forall i,j\in S, j\ne 1, & \forall l: 0\le l\le 2k-1,\\
    \label{sugar}
    p^{(R-1)}_{ij}(z_m)f_j(y_0)&\leq A^R p^{(R-1)}_{b_01}(z_m)f_1(y_0)\nonumber \\
\forall i,j\in     S,&\\
    \label{p}
    {p^{(m+P-1)}_{ij}(y_{2kL})\over p^{(m+P-1)}_{1j}(y_{2kL})}&\leq
q^{-1}\Bigl({K\over \delta}\Bigr)^{m-1}\nonumber \\
 \forall j\in C &\forall i\in S.
\end{align}
If $L=1$, then \eqref{krim} becomes $p_{i\,1}\leq p_{1\,1}$ for all $i\in S$, which is
true by the assumption $p^*_1=p_{1\,1}$ made in the course of constructing the $\mathbf{s}$ sequence (\S\ref{subsec:abs}).
If $L=1$, then \eqref{tsiv} becomes
$${p_{ij}f_j(y_{l+1})\over p_{11}f_1(y_{l+1})}<1-\epsilon \quad \forall i,j\in S, j\ne 1,$$
and thus, since $y_{l+1}\in {\cal X}_1,~0\le l<2k$ in this case, \eqref{tsiv} is true by the definition of
${\cal X}_1$ (\S\ref{subsec:Xl}) (and the fact that $p^*_1=p_{1\,1}$).
Let us next prove \eqref{krim} and \eqref{tsiv} for the case $L>1$.
Consider any $l=0,1,\ldots,2k-1$. Note that the definitions of the $s$-path \eqref{path},
${\cal X}_{s_i}$ (\S\ref{subsec:Xl}), and the fact that $y_{lL+i}\in {\cal X}_{s_i}$ for $1\le i <L$
imply that given observations $y_{Ll+1:L(l+1)-1}$, the path
$s_{1:L-1}$ realizes the maximum in
$p^{(L-1)}_{11}(y_{Ll})$, i.e.
\begin{align}\label{tee}
p_{11}^{(L-1)}(y_{lL})=&p_{1\,s_1}f_{s_1}(y_{lL+1})p_{s_1\,s_2}\cdots \\
                 &\cdots p_{s_{L-2}\,s_{L-1}}f_{s_{L-1}}(y_{(l+1)L-1})p_{s_{L-1}\,1}.\nonumber
\end{align}
(Indeed, $p_{1\,s_1}f_{s_1}(y_{lL+1})p_{s_1\,s_2}\cdots$
\begin{align*}
\cdots
p_{s_{L-2}\,s_{L-1}}f_{s_{L-1}}(y_{(l+1)L-1})p_{s_{L-1}\,1}&=\\
p^*_{s_1}f_{s_1}(y_{lL+1})p^*_{s_2}\cdots
p^*_{s_{L-1}}f_{s_{L-1}}(y_{(l+1)L-1})p^*_{1}&,
\end{align*}
and for $i=1,2,\ldots,L-1$, $p^*_{s_i}f_{s_i}(y_{lL+i})\ge p_{hj}f_{j}(y_{lL+i})$ for any $h,j\in S$.)
Suppose $j\neq 1$ and
$t_{1:L-1}$ realizes $p_{ij}^{(L-1)}(y_{lL}),$ i.e.
\begin{align}\label{kohv}
p_{ij}^{(L-1)}(y_{lL})=&p_{i\,t_1}f_{t_1}(y_{lL+1})p_{t_1\,t_2}\cdots \\ \nonumber
       &\cdots p_{t_{L-2}\,t_{L-1}}f_{t_{L-1}}(y_{(l+1)L-1})p_{t_{L-1}\,j}.
\end{align}
Hence, with $t_0$ and  $t_L$ standing for $i$ and $j$, respectively (and $s_0=s_L=1$),
the left-hand side of \eqref{tsiv} becomes
\begin{align}\label{eq:bigproduct}
&\Bigl({p_{t_0\,t_1}f_{t_1}(y_{lL+1})\over
    p_{s_0\,s_1}f_{s_1}(y_{lL+1})}\Bigr)\Bigl({p_{t_1\,t_2}f_{t_2}(y_{lL+2})\over
    p_{s_1\,s_2}f_{s_2}(y_{lL+2})}\Bigr)\cdots\\
& \Bigl({p_{t_{L-2}\,t_{L-1}}f_{t_{L-1}}(y_{(l+1)L-1})
\over
p_{s_{L-2}\,s_{L-1}}f_{s_{L-1}}(y_{(l+1)L-1})}\Bigr)\Bigl({p_{t_{L-1}\,t_L}f_j(y_{(l+1)L})\over
p_{s_{L-1}\,s_L}f_1(y_{(l+1)L})}\Bigr).\nonumber \end{align}
For $h=1,\ldots,L$ such that  $t_h\ne s_h$,
\begin{eqnarray}\label{kuku}
\lefteqn{{p_{t_{h-1}\,t_h}f_{t_h}(y_{lL+h})\over
p_{s_{h-1}\,s_h}f_{s_h}(y_{lL+h})}<1-\epsilon,~\text{since }y_{lL+h}\in \mathcal{X}_{s_h}.}\hspace{8cm}
\end{eqnarray}
For all other $h$, $s_h=t_h$ and therefore, the left-hand side of \eqref{kuku}
becomes ${p_{t_{h-1}\,t_h}\over p_{s_{h-1}\,s_h}}=\frac{p_{t_{h-1}\,s_h}}{p^*_{s_h}} \leq 1$ (by property \eqref{klass}).
Since the last term of the product \eqref{eq:bigproduct} above does satisfy \eqref{kuku} ($j\ne 1$), \eqref{tsiv} is thus proved.
Suppose next that $t_1,\ldots,t_{L-1}$ realizes $p_{i1}^{(L-1)}(y_{lL})$. With $s_0=1$ and $t_0=i$, similarly to
the previous arguments, we have
$${p_{i\,1}^{(L-1)}(y_{lL})\over p_{1\,1}^{(L-1)}(y_{lL})}=
\prod_{h=1}^{L-1} \Bigl({p_{t_{h-1}\,t_h}f_{t_h}(y_{lL+h})\over
p_{s_{h-1}\,s_h}f_{s_h}(y_{lL+h})}\Bigr){p_{t_{L-1}\,1}\over
p_{s_{L-1}\,1}}\leq 1,$$ implying \eqref{krim}.\\\\
Let us now prove \eqref{sugar}. To that end, note that for all states
$h,i,j\in S$ such that $p_{jh}>0$, it follows from the definitions \eqref{eq:pmaxtrans}
and \eqref{eq:maxtransandratio} that
\begin{equation}\label{A}
{p_{ih}\over p_{jh}}\leq {p^*_h \over p_{jh}}\leq A.\end{equation}
If $R=1$, then \eqref{sugar} becomes $$p_{ij}f_j(y_0)\leq
Ap_{b_01}f_1(y_0).$$ By the definition of ${\cal X}_1$ (recall that
$y_0\in {\cal X}_1$), we have that for every $i,j\in S$
$p_{ij}f_j(y_0)\le p^*_1f_1(y_0)$. Using \eqref{A} with $h=1$ and $j=b_{0}$, we get
$p^*_1f_1(y_0)\leq Ap_{b_0\,1}f_1(y_0)$ ($p_{b_0\,1}>0$ by the construction of $\mathbf{b}$
\S\ref{subsec:abs}). Putting these all together, we obtain
\begin{equation*}
p_{ij}f_j(y_0)<p^*_1f_1(y_0)\leq Ap_{b_0 1}f_1(y_0),~\text{as required.}
\end{equation*}
Consider the case $R>1$. Let $t_{1:R-1}$ be a path that
realizes $p^{(R-1)}_{ij}(z_m)$, i.e. $p^{(R-1)}_{ij}(z_m)=$
$$p_{i\,t_1}f_{t_1}(y'_1)p_{t_1\, t_2}f_{t_2}(y'_2)\cdots
p_{t_{R-2}\,t_{R-1}}f_{t_{R-1}}(y'_{R-1})p_{t_{R-1}j}.$$
By 
the definition of $\mathcal{X}_l$ (\S\ref{subsec:Xl}) and
the facts that $y'_r\in \mathcal{X}_{b_r}$, $r=1,2,\ldots,R-1$, and $y_0\in\mathcal{X}_1$,
we have
\begin{multline}\label{eq:claim3proof}
p^{(R-1)}_{ij}(z_m)f_j(y_0)\le p^*_{b_1}f_{b_1}(y'_1)p^*_{b_2}f_{b_2}(y'_2)\cdots\\
p^*_{b_{R-1}}f_{b_{R-1}}(y'_{R-1})p^*_1f_1(y_0).\end{multline}
Now, by the construction of $\mathbf{b}$ (\S\ref{subsec:abs}), $p_{b_{r-1}\,b_r}>0$ for $r=1,\ldots,R$, ($b_R=1$).
Thus,  the argument behind \eqref{A} applies here to bound the right-hand side of \eqref{eq:claim3proof} from above
by
\begin{multline*}
Ap_{b_0\,b_1}f_{b_1}(y'_1)Ap_{b_1\,b_2}f_{b_2}(y'_2)\cdots \\
Ap_{b_{R-2}\,b_{R-1}}f_{b_{R-1}}(y'_{R-1})Ap_{b_{R-1}\,1}f_1(y_0)=\\ A^Rp^{(R-1)}_{b_0\,1}(z_m)f_1(y_0),
~\text{as required.}
\end{multline*}
Let us now prove \eqref{p}.
If $m=1$ then \eqref{p} becomes
\begin{equation}\label{teem1}
p^{(P)}_{ij}(y_{2kL})\leq p^{(P)}_{1j}(y_{2kL})q^{-1}\quad \forall
j\in C \forall i\in S.
\end{equation}
If $P=0$, then \eqref{teem1} reduces to $p_{ij}\leq p_{1j}q^{-1}$ which is true, because in this case
the state $q_1=q_T=1$ belongs to $C$ (\S\ref{subsec:abs}) and
$p_{1j}q^{-1}\geq 1$ (\eqref{maxprobpath}, \eqref{eq:qpositive} with $m=1$).
To see why \eqref{teem1} is true with $P\geq 1$, note  
that by the same argument as used for proving \eqref{krim} and \eqref{tsiv}, we now get
$ \forall h,l\in S$
\begin{equation}\label{tele}
p_{1\,a_P}^{(P-1)}(y_{2kL})f_{a_P}(y^{''}_P)\geq
p_{h',l}^{(P-1)}(y_{2kL})f_{l}(y^{''}_P).\end{equation}
Also, since $a_P=q_1\in C$ (\S\ref{subsec:abs}),
$p_{a_P\,j}q^{-1}\geq 1$ (\eqref{maxprobpath}, \eqref{eq:qpositive} with $m=1$). Thus
$p_{i\,j}^{(P)}(y_{2kL})=$
\begin{align*}
&\stackrel{\mathrm{by~\eqref{eq:prrecurse}}}{=}
\max_{l\in S} p_{i\,l}^{(P-1)}(y_{2kL})f_l(y^{''}_P)p_{l\,j}\\
&\stackrel{\mathrm{by~\eqref{tele}}}{\le} p^{(P-1)}_{1a_P}(y_{2kL})f_{a_P}(y''_P)\max_{l\in S}p_{l\,j}\\
&  \le p^{(P-1)}_{1\,a_P}(y_{2kL})f_{a_P}(y''_P)\\
&\le
p_{1\,a_P}^{(P-1)}(y_{2kL})f_{a_P}(y^{''}_P)p_{a_P\,j}q^{-1}\stackrel{\mathrm{by~\eqref{eq:prrecurse}}}{\le}
p^{(P)}_{1\,j}(y_{2kL})q^{-1}.
\end{align*}
For $m>1$, let $t_{1:m-1}$ be a path realizing $p^{(m-1)}_{h\,j}(y^{''}_P)$.
Thus, $p^{(m-1)}_{h\,j}(y^{''}_P)=$
\begin{align}\label{eq:auxbound1}
&=p_{h\,t_1}f_{t_1}(z'_1)p_{t_1\,t_2}f_{t_2}(z'_2)\cdots
f_{t_{m-1}}(z'_{m-1})p_{t_{m-1}j}\nonumber \\
&< K^{m-1}.
\end{align}
(This is true since $z'_r\in \mathcal{Z}$ for $r=1,2,\ldots,m-1$ (\S\ref{subsec:Z})
and thus, for $p^{(m-1)}_{h\,j}(y^{''}_P)$ to be positive it is
necessary that  $t_r\in C$, $r=1,\ldots,m-1$, implying $f_{t_r}(z'_r)<K$.) Now,
let $t_{1:m-1}$ realize $p^{(m-1)}_{a_P\,j}(y^{''}_P)$,
which is clearly positive, with $t_r\in C$, $r=1,\ldots,m-1$
($z'_r\in \mathcal{Z}$ for $r=1,2,\ldots,m-1$), and $a_P, j\in C$
(recall the positivity assumption on $\mathbb{Q}^m$, \S\ref{subsec:k}).
We thus have $p^{(m-1)}_{a_P\,j}(y^{''}_P)=p_{a_P\,t_1}f_{t_1}(z'_1)p_{t_1\,t_2}f_{t_2}(z'_2)\cdots
f_{t_{m-1}}(z'_{m-1})p_{t_{m-1}j}\ge $
\begin{eqnarray}\label{eq:auxbound2}
\lefteqn{\ge q^*_{a_P\,j}f_{t_1}(z'_1)f_{t_2}(z'_2)\cdots
f_{t_{m-1}}(z'_{m-1})>q \delta^{m-1}.}\hspace{8cm}
\end{eqnarray}
Combining the  bounds of \eqref{eq:auxbound1} and \eqref{eq:auxbound2} ($q>0$, \eqref{eq:qpositive}), we obtain:
\begin{equation} \label{eq:teem1further}
p^{(m-1)}_{h\,j}(y^{''}_P) < p^{(m-1)}_{a_P\,j}(y^{''}_P)\Bigl({K\over \delta}\Bigr)^{m-1}/q.
\end{equation}
Finally, $p_{ij}^{(P+m-1)}(y_{2kL})=$
\begin{align*}
&\stackrel{\mathrm{by~\eqref{eq:prrecurse}}}{=}\max_{l\in S}
p_{il}^{(P-1)}(y_{2kL})f_l(y^{''}_P)p^{(m-1)}_{lj}(y^{''}_P)\\
&\stackrel{\mathrm{by~\eqref{tele},~\eqref{eq:teem1further}}}{<}
p_{1\,a_P}^{(P-1)}(y_{2kL})f_{a_P}(y^{''}_P)p^{(m-1)}_{a_P\,j}(y^{''}_P)\left({K\over
\delta}\right)^{m-1}\hspace*{-5mm}/q\\ &\stackrel{\mathrm{by~\eqref{eq:prrecurse}}}{\le}
 p_{1j}^{(P+m-1)}(y_{2kL})\left({K\over
\delta}\right)^{m-1}/q.
\end{align*}
\subsubsection{$\gamma_j\le const \times \gamma_1$}
\label{subsec:gammabound}
Combining \eqref{kreeka2}, 
\eqref{betas}, 
and \eqref{sugar}, we see that
for every state $j\in S$,
\begin{align*}
\gamma_j&\stackrel{\mathrm{by~\eqref{kreeka2}}}{=}\beta_{i_\gamma(j)}p_{i_\gamma(j)\,j}^{(R-1)}(z_m)f_j(y_0)\\
&\stackrel{\mathrm{by~\eqref{sugar}}}{\leq} \beta_{i_\gamma(j)}
p_{b_0\,1}^{(R-1)}(z_m)f_1(y_0)A^R\\ & \stackrel{\mathrm{by~\eqref{betas}}}{\leq}
 q^{-1}\Bigl({K\over
\delta}\Bigr)^{m}A^R\beta_{b_0}p_{b_0\,1}^{(R-1)}(z_m)f_1(y_0)\\ & \leq U
\max_{i\in S} \beta_i p_{i\,1}^{(R-1)}(z_m)f_1(y_0)\stackrel{\mathrm{by~\eqref{kreeka2}}}{=}U\gamma_1,
\end{align*}
where
\begin{equation}\label{eq:defU}
U\stackrel{\mathrm{def}}{=}q^{-1}\Bigl({K\over
\delta}\Bigr)^{m}A^R.\end{equation}
Hence
\begin{equation}\label{con}
\gamma_j\leq U\gamma_1 \quad  \forall j\in S. \end{equation}
\subsubsection{Further bounds on likelihoods}
\label{subsec:gammaboundmore}
Let $l\ge 0$ and $n>0$ be integers such that  $l+n \leq 2k$ but arbitrary otherwise.
Expanding $p_{1\,1}^{(nL-1)}(y_{lL})$ recursively according with \eqref{eq:prrecurse}, we obtain
\begin{multline}\label{prekohus}
p_{1\,1}^{(nL-1)}(y_{lL})=\max_{i_{1:n-1}\in S^{n-1}} p_{1\,i_1}^{(L-1)}(y_{lL})f_{i_1}(y_{(l+1)L}) \times \\
\times p_{i_1\,i_2}^{(L-1)}(y_{(l+1)L})f_{i_2}(y_{(l+2)L}) \cdots
p_{i_{n-2}\,i_{n-1}}^{(L-1)}(y_{(l+n-2)L}) \times \\
\shoveright{\times f_{i_{n-1}}(y_{(l+n-1)L})p_{i_{n-1}\,1}^{(L-1)}(y_{(l+n-1)L}).}
\end{multline}
Since for any $i_1\in S$, $p_{1\,i_1}^{(L-1)}(y_{lL})f_{i_1}(y_{(l+1)L})\le p_{1\,1}^{(L-1)}(y_{lL})f_{1}(y_{(l+1)L})$,
as well as
\begin{multline*}
p_{i_{r-1}\,i_r}^{(L-1)}(y_{(l+r-1)L})f_{i_r}(y_{(l+r)L})\stackrel{\mathrm{by~\eqref{tsiv}}}{\leq}  \\
p_{1\,1}^{(L-1)}(y_{(l+r-1)L})f_1(y_{(l+r)L}),~ r=2,\ldots,n-1,
\end{multline*}
and since for any $i_{n-1}\in S$
\begin{equation*}p_{i_{n-1}\,1}^{(L-1)}(y_{(l+n-1)L})\stackrel{\mathrm{by~\eqref{krim}}}{\leq}
 p_{1\,1}^{(L-1)}(y_{(l+n-1)L}),
\end{equation*}
maximization \eqref{prekohus} above is achieved as
\begin{align}\label{kohus}
&\text{follows:}~p_{1\,1}^{(nL-1)}(y_{lL})=\\
&p_{1\,1}^{(L-1)}(y_{lL})f_1(y_{(l+1)L})p_{11}^{(L-1)}(y_{(l+1)L})f_1
(y_{(l+2)L}) \cdots \nonumber \\
& \cdots
p_{1\,1}^{(L-1)}(y_{(l+n-2)L})f_1(y_{(l+n-1)L})p_{1\,1}^{(L-1)}(y_{(l+n-1)L}).\nonumber
\end{align}
Now, we replace state $1$ by generic states $i,j\in S$ on the both ends of the paths in \eqref{prekohus}
and repeat the above arguments. Thus, also using \eqref{kohus}, we arrive at bound \eqref{liim} below:
\begin{align}\label{liim}
p_{ij}^{(nL-1)}(y_{lL})f_j(y_{(l+n)L})&\leq \nonumber \\
\prod_{u=l+1}^{l+n}
p_{11}^{(L-1)}(y_{(u-1)L})f_1(y_{uL})&\stackrel{\mathrm{by~\eqref{kohus}}}{=}\nonumber  \\
p_{11}^{(nL-1)}(y_{lL})f_1(y_{(l+n)L}) \quad \forall i,j\in S.
\end{align}

In particular, \eqref{liim} states $\forall i,j\in S$
\begin{equation}\label{pulk}
p_{ij}^{(kL-1)}(y_0)f_j(y_{kL})\leq
p_{11}^{(kL-1)}(y_0)f_1(y_{kL}).
\end{equation}
\subsubsection{$\eta_j\le const \times \eta_1$}
\label{subsec:etabound}
In order to see
\begin{equation}\label{gamma}
\eta_j\leq U\eta_1 \quad \forall j\in S,
\end{equation}
\begin{align*}
\text{note:}~\eta_j&\stackrel{\mathrm{\eqref{kreeka3}}}{=}\max_{i\in S} \gamma_i p_{i\,j}^{(kL-1)}(y_0)f_j(y_{kL})\\
&\stackrel{\mathrm{by~\eqref{pulk}}}{\leq}
\max_{i\in S} \gamma_i p_{1\,1}^{(kL-1)}(y_0)f_1(y_{kL})\stackrel{\mathrm{by~\eqref{con}}}{\leq}\\
   & \stackrel{\mathrm{by~\eqref{con}}}{\leq}U\gamma_1 p_{1\,1}^{(kL-1)}(y_0)f_1(y_{kL})
 \stackrel{\mathrm{by~\eqref{kreeka3}}}{\le } U\eta_1.
\end{align*}
\subsubsection{A representation of  $\eta_1$}
\label{subsec:etaone}
Recall that $k$, the number of cycles in the $s$-path, was chosen sufficiently large for
\eqref{eq:largek} to hold (in particular, $k>1$).
We now prove that there exists $\kappa\in \{1,\ldots,k-1\}$ such that
\begin{equation}\label{on}
\eta_1=\delta_{1}(y_{\kappa L})p_{1\,1}^{((k-\kappa)L-1)}(y_{\kappa L})f_1(y_{kL}).
\end{equation}
The relation \eqref{on} states that (given observations $x_{1:u}$) a maximum-likelihood path
(from time $1$, observation $x_1$) to time $u-m-P-kL$ (observation $y_{kL}$)
goes through state $1$ at time $u-m-P-2kL+\kappa L$, that is when $y_{\kappa L}$ is observed.

To see this, suppose no such $\kappa$ existed.
Then,
applying \eqref{eq:prrecurse} to \eqref{kreeka3} and recalling that
$\delta_{1}(y_{\kappa L})$ is introduced in \eqref{eq:shortscores}, we would have
\begin{multline*}
\eta_1=\gamma_{j_\eta(1)}p_{j_\eta(1)\,j_1}^{(L-1)}(y_0)f_{j_1}(y_L)
p_{j_1\,j_2}^{(L-1)}(y_L)\times \\
\times f_{j_2}(y_{2L})p_{j_2\,j_3}^{(L-1)}(y_{2L})\cdots
p_{j_{k-1}\,1}^{(L-1)}(y_{(k-1)L})f_1(y_{kL})
\end{multline*}
for some $j_1\ne 1,\ldots,j_{k-1}\ne 1$. Furthermore, this would imply 
$\eta_1<$
\begin{align}\label{jalle}
\nonumber
&\stackrel{\mathrm{by~\eqref{tsiv},~\eqref{krim}}}{<}\gamma_{j_\eta(1)}(1-\epsilon)^{k-1}
\prod_{i=1}^{k}p_{1\,1}^{(L-1)}(y_{(i-1)L})f_1(y_{iL})\\
\nonumber
& \hspace*{-2mm}\stackrel{\mathrm{by~\eqref{eq:largek}}}{<}
\gamma_{j_\eta(1)}q^2\left({\delta\over K}\right)^{2m}A^{-R} \prod_{i=1}^{k}p_{1\,1}^{(L-1)}(y_{(i-1)L})f_1(y_{iL})
\\
\nonumber
& \stackrel{\mathrm{by~\eqref{con}}}{\le}
\gamma_1Uq^2\left ({\delta\over K}\right )^{2m}A^{-R} \prod_{i=1}^{k}p_{1\,1}^{(L-1)}(y_{(i-1)L})f_1(y_{iL})
\\ & \stackrel{\mathrm{by~\eqref{eq:defU}}}{=}
\gamma_1q\left ({\delta\over K}\right )^m \prod_{i=1}^{k}p_{1\,1}^{(L-1)}(y_{(i-1)L})f_1(y_{iL})\nonumber \\
& <\gamma_1 \prod_{i=1}^{k}p_{1\,1}^{(L-1)}(y_{(i-1)L})f_1(y_{iL}).
\end{align}
(The last inequality follows from $q\le 1$ \eqref{eq:qpositive} and $\delta<K$, \S\ref{subsec:Z}.)
On the other hand, by definition \eqref{kreeka3} (and $k-1$-fold application of \eqref{eq:prrecurse}),
$\eta_1\ge \gamma_1 \prod_{i=1}^{k}p_{1\,1}^{(L-1)}(y_{(i-1)L})f_1(y_{iL}),$ which evidently contradicts
\eqref{jalle} above.
Therefore, $\kappa$ satisfying \eqref{on} and $1\le \kappa< k$, does exist. 
\subsubsection{An implication  of  \eqref{kohus} and \eqref{on} for $\delta_1(y_{lL})$}
\label{subsec:deltaone}
Clearly, the arguments of the previous section (\S\ref{subsec:etaone}) are valid
if $k$ is replaced by any $l\in\{k,\ldots,2k\}$. Hence the following generalization of
\eqref{on}: For some $\kappa(l)<l$
\begin{equation}\label{preojakaar}
\hspace*{-3mm}\delta_1(y_{lL})=\delta_1(y_{\kappa(l)L})p^{((l-\kappa(l))L-1)}_{11}(y_{\kappa(l)L})f_1(y_{lL}).
\end{equation}
We apply \eqref{preojakaar} recursively, starting with  $\kappa^{(0)}\stackrel{\mathrm{def}}{=}l$ and returning
$\kappa^{(1)}\stackrel{\mathrm{def}}{=}\kappa(l)<l$. If $\kappa^{(1)}\leq k$, we stop, otherwise we substitute $\kappa^{(1)}$
for $l$, and obtain $\kappa^{(2)}\stackrel{\mathrm{def}}{=}\kappa(l)<\kappa^{(1)}$, and so, on until $\kappa^{(j)}\le k$ for some $j>0$.
Thus, $\delta_1(y_{lL})=$
\begin{multline}\label{preojakaar2}\hspace*{-4mm}
=\delta_1(y_{\kappa^{(j)}L})p^{((\kappa^{(j-1)}-\kappa^{(j)})L-1)}_{11}(y_{\kappa^{(j)}L})
f_1(y_{\kappa^{(j-1)}L})\cdots\\ p^{((l-\kappa^{(1)})L-1)}_{11}(y_{\kappa^{(1)}L})
f_1(y_{lL}).
\end{multline}
Applying \eqref{kohus} to the appropriate factors of the right-hand side of \eqref{preojakaar2} above,
we obtain:
\begin{multline}\label{preojakaar3}
\delta_1(y_{lL})  = \delta_1(y_{\kappa^{(j)}L})p^{(L-1)}_{11}(y_{\kappa^{(j)}L})f_1(y_{(\kappa^{(j)}+1)L})\cdots \\
       p^{(L-1)}_{11}(y_{(k-1)L})f_1(y_{kL})\cdots
  p^{(L-1)}_{11}(y_{kL})f_1(y_{(k+1)L})\cdots\\
            p^{(L-1)}_{11}(y_{(\kappa^{(j-1)}-1)L})f_1(y_{\kappa^{(j-1)}L})\cdots \\
   p^{(L-1)}_{11}(y_{(\kappa^{(1)}-1)L})f_1(y_{\kappa^{(1)}L})\cdots \\p^{(L-1)}_{11}(y_{(l-1)L})f_1(y_{lL}).
\end{multline}
Also, according to \eqref{kohus}, \begin{multline*}
 \delta_1(y_{\kappa^{(j)}L})p^{(L-1)}_{11}(y_{\kappa^{(j)}L})f_1(y_{(\kappa^{(j)}+1)L})\cdots\\
p^{(L-1)}_{11}(y_{(k-1)L}) =\delta_1(y_{\kappa^{(j)}L})p^{((k-\kappa^{(j)})L-1)}_{11}(y_{\kappa^{(j)}L}).
\end{multline*}
At the same time,
\begin{equation}\label{preojakaar4}
\hspace*{-3mm}\delta_1(y_{\kappa^{(j)}L})p^{((k-\kappa^{(j)})L-1)}_{11}(y_{\kappa^{(j)}L})f_1(y_{kL})
\stackrel{\mathrm{by~\eqref{eq:prrecurse}}} \le \eta_1.
\end{equation}
However, we cannot have the strict inequality in \eqref{preojakaar4} above since that, by virtue of
\eqref{preojakaar3},  would contradict maximality of $\delta_1(y_{lL})$.
We have thus arrived at $\delta_1(y_{lL})=\eta_1p^{(L-1)}_{11}(y_{{k}L})f_1(y_{(k +1)L})\cdots $
\begin{equation}\label{ojakaar}
 \cdots p^{(L-1)}_{11}(y_{(l-1)L})f_1(y_{lL}).
\end{equation}

In summary, for any $l\ge k$ and $l\le 2k$ there exists a realization of  $\delta_1(y_{lL})$ that goes
through state $1$ every time when $y_{iL}$, $i=k,\ldots, l$, is observed.
\subsubsection{$y_{kL}$ is a $(kL+m+P)$-order 1-node}
\label{subsec:onenode}
In \S\ref{subsec:auxproof35}, we will prove that
for any  $i\in S,i\ne 1,$ and any$ j\in C$,
\begin{equation}\label{lill}
\eta_ip_{ij}^{(kL+m+P-1)}(y_{kL})\leq
\eta_1p_{1j}^{(kL+m+P-1)}(y_{kL}),
\end{equation}
which implies that $y_{kL}$ is a  1-node of order $kL+m+P$. Indeed,
let $l\in S$ be arbitrary. Since $f_j(z'_m)=0$ for every $j\in S\setminus C$, any maximum
likelihood path to state $l$ at time $u+1$ (observation $x_{u+1}$) must go
through a state in $C$ at time $u$ (observation $x_u=z'_m$.) Formally,
\begin{eqnarray*}
\lefteqn{\eta_ip_{il}^{(kL+m+P)}(y_{kL})=}\\
&=& \max_{j\in
S}\eta_ip_{ij}^{(kL+m+P-1)}(y_{kL})f_j(z'_m)p_{jl}\\
&=& \max_{j\in C}\eta_ip_{ij}^{(kL+m+P-1)}(y_{kL})f_j(z'_m)p_{jl}    \\
&\stackrel{\mathrm{by~\eqref{lill}}}{\le}&  \max_{j\in C}\eta_1p_{1j}^{(kL+m+P-1)}(y_{kL})f_j(z'_m)p_{jl}\\
&\stackrel{\mathrm{by~\eqref{eq:prrecurse}}}{=}& \eta_1p_{1l}^{(kL+m+P)}(y_{kL}).
\end{eqnarray*}
Therefore, by Definition~\ref{rnode} $y_{kL}$ is a 1-node of order $kL+m+P$.
\subsubsection{Proof of  \eqref{lill}}
\label{subsec:auxproof35}
Let $i\in S$ and $j\in C$ be arbitrary. Let state $j^*\in S$ be such that
$p_{i\,j}^{(kL+m+P-1)}(y_{kL})=$
$$p_{i\,j^*}^{(kL-1)}(y_{kL})f_{j^*}(y_{2kL})p^{(m+P-1)}_{j^*\,j}(y_{2kL})=$$
$\nu(i,j^*)p^{(m+P-1)}_{j^*\,j}(y_{2kL}),$ where
$$\nu(i,j)\stackrel{\mathrm{def}}{=}p_{ij}^{(kL-1)}(y_{kL})f_j(y_{2kL}),\quad\text{for all}~i,j\in S.$$
We consider the following two cases separately:
\begin{enumerate}[1.]
\item 
There exists a path realizing $p_{i\,j^*}^{(kL-1)}(y_{kL})$ and going
through state 1 at the time of observing $y_{lL}$ for some
$l \in \{k,\ldots, 2k\}$.
$p_{i\,j^*}^{(kL-1)}(y_{kL})=$
\begin{equation}\label{ikka}
p^{((l-k) L-1)}_{i\,1}(y_{kL})f_1(y_{lL})p^{((2k-l)L-1)}_{1\,j^*}(y_{lL}).
\end{equation}
Equation \eqref{ikka} above together with the fundamental recursion \eqref{eq:prrecurse} yields
the following:
\begin{eqnarray}
    \label{eq:ikka2}\nonumber
\lefteqn{\eta_ip_{i\,j^*}^{(kL-1)}(y_{kL})=}\\
&& \stackrel{\mathrm{by~\eqref{ikka}}}{=}
\eta_ip^{((l-k) L-1)}_{i\,1}(y_{kL})f_1(y_{lL})p^{((2k-l)L-1)}_{1\,j^*}(y_{lL})\nonumber
 \\
&& \stackrel{\mathrm{by~\eqref{eq:shortscores},~\eqref{eq:prrecurse}}}{\le}
\delta_1(y_{lL})p^{((2k-l)L-1)}_{1\,j}(y_{lL}).
  \end{eqnarray}
At the same time, the right hand-side of \eqref{eq:ikka2} can be expressed as follows:
\begin{eqnarray}
    \label{eq:ikka3}\nonumber
\lefteqn{\delta_1(y_{lL})p^{((2k-l)L-1)}_{1\,j^*}(y_{lL})}\\
&&\stackrel{\mathrm{by~\eqref{ojakaar}}}{=}
\eta_1p_{1\,1}^{((l-k)L-1)}(y_{kL})f_1(y_{lL})p^{((2k-l)L-1)}_{1\,j^*}
\nonumber\\
&&\stackrel{\mathrm{by~\eqref{kohus}}}{=}\eta_1p_{1\,j^*}^{(kL-1)}(y_{kL}).
  \end{eqnarray}
Therefore, if there exists $l \in \{k,\ldots, 2k\}$ such that \eqref{ikka} holds,
we have by virtue of \eqref{eq:ikka2} and \eqref{eq:ikka3}:
$\eta_ip_{i\,j^*}^{(kL-1)}(y_{kL})\leq\eta_1p_{1\,j^*}^{(kL-1)}(y_{kL}),$ that is
\begin{equation}\label{ometigi}
\eta_i\nu(i,j^*)\leq \eta_1\nu(1,j^*).
\end{equation}
\begin{eqnarray*}
\lefteqn{\hspace*{-1cm}\text{Hence,}~\eta_ip_{i\,j}^{(kL+m+P-1)}(y_{kL})=}\\
&\stackrel{\mathrm{by~\eqref{ikka}}}{=}&
\eta_i\nu(i,j^*)p_{j^*\,l}^{(m+P-1)}(y_{2kL})\\
&\stackrel{\mathrm{by~\eqref{ometigi}}}{\le}&
\eta_1\nu(1,j^*)p_{j^*\,j}^{(m+P-1)}(y_{2kL})\\
&\stackrel{\mathrm{by~\eqref{eq:prrecurse}}}{\le}&
\eta_1p_{1\,j}^{(kL+m+P-1)}(y_{kL})
\end{eqnarray*}
and \eqref{lill} holds.
\item 
Assume now that no path exists to satisfy \eqref{ikka}.
Argue as for \eqref{jalle} to obtain $\nu(i,j^*)<$
\begin{equation}\label{jalle2}
(1-\epsilon)^{k-1} \prod_{n=k+1}^{2k}p_{1\,1}^{(L-1)}(y_{(n-1)L})f_1(y_{nL}).
\end{equation} By \ref{kohus}, the (partial likelihood) product in the
right-hand side of \eqref{jalle2} equals $\nu(1,1)$. Thus,
\begin{eqnarray}\label{lopp}
\lefteqn{\eta_i\nu(i,j^*) p^{(m+P-1)}_{j^*\,j}(y_{2kL})<}\\ \nonumber
&&\stackrel{\mathrm{by~\eqref{jalle2}}}{<}
\eta_{i}(1-\epsilon)^{k-1}\nu(1,1)p^{(m+P-1)}_{j^*\,j}(y_{2kL})
 \\ \nonumber
&& \stackrel{\mathrm{by~\eqref{eq:largek}}}{<} \eta_{i} q^2\left({\delta\over K}\right)^{2m}A^{-R}
\nu(1,1)p^{(m+P-1)}_{j^*\,j}(y_{2kL}) \\ \nonumber
&& \stackrel{\mathrm{by~\eqref{eq:defU},~\eqref{gamma}}}{\le}\eta_1 q\left({\delta\over K}\right)^{m}
\nu(1,1)p^{(m+P-1)}_{j^*\,j}(y_{2kL}).
\end{eqnarray}
Hence, for every $j'\in S$,
\begin{eqnarray*}
\lefteqn{\eta_i\nu(i,j')p^{(m+P-1)}_{j'\,j}(y_{2kL})}\\
&&\stackrel{\mathrm{by~\eqref{ikka}}}{\le}
\eta_i\nu(i,j^*)p^{(m+P-1)}_{j^*\,j}(y_{2kL})\stackrel{\mathrm{by~\eqref{lopp}}}{<}\\
&&\stackrel{\mathrm{by~\eqref{lopp}}}{<}
\eta_1 q\left({\delta\over K}\right)^{m}\nu(1,1)p^{(m+P-1)}_{j^*\,j}(y_{2kL})\\
&&\stackrel{\mathrm{by~\eqref{p}}}{\le}
\eta_1 \left({\delta\over K}\right)\nu(1,1)p^{(m+P-1)}_{1\,j}(y_{2kL})\\
&&<
\eta_1 \nu(1,1)p^{(m+P-1)}_{1\,j}(y_{2kL})\\
&&\stackrel{\mathrm{by~\eqref{eq:prrecurse}}}{\le} \eta_1p_{1\,j}^{(kL+m+P-1)}(y_{kL}),
\end{eqnarray*}
\noindent which, by virtue of \eqref{eq:prrecurse}, implies \eqref{lill}.
\end{enumerate}
\subsubsection{Completion of the $s$-path to $q_{1:M}$ and conclusion}
\label{subsec:M}
Finally, let $$M=2m+2Lk+P+R+2,\quad r=kL+P+m,~ l=1.$$
Recall from \S\ref{subsec:abs} that
$b_0\in C$. Since all the entries of $\mathbb{Q}^m$ are positive, there exists a path
$v_{0:m-1},b_0\in C$ such that $p_{v_i\,v_{i+1}}>0$ and
$p_{v_{m-1}b_0}>0$.
Similarly, there must exist a path $u_{1:m}\in C$ such that $p_{u_i\,u_{i+1}}>0$
$\forall i=1,\ldots, m-1$ and $p_{a_P\,u_1}>0$ (recall that $a_P\in C$).
Hence, by these, and the constructions of \S\ref{subsec:spath},
all of the  transitions of the following sequence occur with positive probabilities.
\begin{equation}\label{pathbig}
q_{1:M}\stackrel{\mathrm{def}}{=}(v_{0:m-1},b_{0:R},s_{1:2Lk},a_{1:P},u_{1:m}).\end{equation}
Clearly, the actual probability of observing $q_{1:M}$ is positive, as required.
By the constructions of \S\S\ref{subsec:Xl}-\ref{subsec:abs}, the conditional
probability of $B$ below, given $q_{1:M}$, is evidently positive, as required.
\begin{multline*}
B\stackrel{\mathrm{def}}{=}{\cal Z}^{m+1}\times {\cal X}_{b_1}\times \cdots \times {\cal X}_{b_{R-1}}\times {\cal X}_1\times {\cal X}_{s_1} \times \\
\cdots \times {\cal X}_{s_{2kL-1}}\times {\cal X}_1\times {\cal
X}_{a_1}\times\cdots \times {\cal X}_{a_{P}}\times {\cal Z}^{m}.
\end{multline*}
Finally, since the sequence \eqref{blokk} below was chosen from $B$
arbitrarily (\S\ref{subsec:barrier})
and has been shown to be an $l$-barrier of order $r$, this completes the proof of the Lemma.
\begin{equation*}
(z_{0:m},y'_{1:R-1},y_{0:2Lk},y^{''}_{1:P},z'_{1:m})\in B.
\hspace*{1.5cm}\eqref{blokk}
\end{equation*}
\subsection{Proof of Lemma \ref{separated}}\label{sec:separatedproof}
\begin{proof}
We use the notation of the previous proof in \S\ref{sec:proofneljas} and
consider the following two distinct situations: First (\S\ref{subsec:donedeal}), all barriers from $B$
as constructed in  the proof of Lemma \ref{separated} are already separated. Obviously,
there is nothing to do in this case.  The second situation (\S\ref{subsec:needmore}) is complementary,
in which case a simple extension will immediately ensure separation.

\subsubsection{All $y\in B$ are already separated}
\label{subsec:donedeal}
Recall the definition of ${\cal Z}$ from \S\ref{subsec:Z}.
Consider the two cases in the definition separately.
First, suppose ${\cal Z}=\hat{{\cal Z}}\backslash (\cup_{l\in S}{\cal X}_l)$,
in which case ${\cal Z}$ and ${\cal X}_l$ are disjoint for every $l\in S$.
This implies that 
every barrier \eqref{blokk} 
is already separated. Indeed, for any $w$, $1\le w\le r$,
and for any $y\in B$, the fact that $y_{M-\max(m,w)}\not\in {\cal Z}$, for example,
makes it impossible for $(y'_{1:w},y_{1:M-w})\in B$ for any $y'_{1:w}\in \mathcal{X}^w$.
Consider now the case when ${\cal Z}=\hat{{\cal Z}}\cap {\cal X}_s$
for some $s\in C$. Then
\begin{multline}
B\subset {\cal X}_s^{m+1}\times
{\cal X}_{b_1}\times \cdots \times {\cal X}_{b_{R-1}}\times {\cal
X}_1\times {\cal X}_{s_1} \times  \cdots \\  {\cal
X}_{s_{2kL-1}}\times {\cal X}_1\times {\cal X}_{a_1}\times\cdots
\times {\cal X}_{a_{P-1}}\times {\cal X}_s^{m+1}.
\end{multline}
Let $y\in B$ be arbitrary.
Assume first  $L>1$. By construction (\S\ref{subsec:abs}), the states $s_1,\ldots,s_L$ are all distinct.
We now show that $(y'_{1:w},y_{1:M-w})\not \in B$ for any $y'_{1:w}\in \mathcal{X}^w$
when $1\le w\le r$.  Note that the sequence
$$q_{m+2: m+R+2kL+P+1}=(b_{1:R-1},1,s_{1:2kL-1},1,a_{1:P-1},s)$$
is such that no two consecutive states are equal.  It is  straightforward to
verify that there exist indices $j$, $0\le j\le m-1$, such that, when shifted $w$ positions to the right,
the pair $y_{j+1\,j+2}\in \mathcal{X}_s^2$ would at the same time have to belong to
$\mathcal{X}_{q_{j+1+w}}\times \mathcal{X}_{q_{j+2+w}}$ with $m+1\le j+1+w<j+2+w\le m+R+2kL+1+P$. This
is clearly a contradiction since  $\mathcal{X}_{q_{j+1+w}}$ and  $\mathcal{X}_{q_{j+2+w}}$ are disjoint
for that range of indices $j$. A verification of the above fact simply amounts to verifying that
the inequality $\max(0,m-w)\le j\le \min(m-1,m+R+2kL-1+P-w)$ is consistent for any $w$ from the admissible
range:
\begin{enumerate}[i.)]
\item When $0\ge m-w$, $m-1\le m+R+2kL-1+P-w$ ($m\le w \le \min(r,R+2kL+P)$), $0\le j \le m-1$
is evidently consistent.
\item When $0\ge m-w$, $m-1> m+R+2kL-1+P-w$ ($\max(m,R+2kL+P)\le w \le r$), $0\le j \le m+R+2kL-1+P-w$
is also consistent since $m+R+2kL-1+P-r=R+kL-1\ge 0$.
\item When $0< m-w$, $m-1\le m+R+2kL-1+P-w$ ($1\le w \le \min(m-1,R+2kL+P)$), $m-w\le j \le m-1$ is
consistent since $w\ge 1$.
\item When $0< m-w$, $m-1> m+R+2kL-1+P-w$ ($\max(1,R+2kL+P-1)\le w < m$), $m-w\le j \le m+R+2kL-1+P-w$ is
consistent since $R+2kL-1\ge 0$.
\end{enumerate}

Next consider the case of $L=1$ but $s\neq 1$ (that is, $P>0$). Then
$B\subset {\cal X}_s^{m+1}\times{\cal X}_{b_1}\times \cdots $
$$\times {\cal X}_{b_{R-1}}\times {\cal
X}_1^{2k+1}\times {\cal X}_{a_{1}}\times \cdots \times {\cal
X}_{a_{P-1}}\times {\cal X}_s^{m+1}.$$   If $s\ne 1$, then also $b_i\ne 1$,
$i=1,\ldots, R-1$ and $a_i\ne 1$, $i=1,\ldots, P-1$. To
see that $y$ is separated in this case, simply note
that $y_{M-max(w,m+1)}\not\in\mathcal{X}_s$ for any admissible $w$.

\subsubsection{Barriers $y\in B$ need not be separated}
\label{subsec:needmore}
Finally, we consider the case when $L=1$ and $s=1$ (where $s\in C$ is such that
$\mathcal{Z}=\hat{\mathcal{Z}}\cap \mathcal{X}_s$).
This implies that $P=0$, $1\in C$, and $p_{1\,1}>0$, which in turn implies that $R=1$, and
$$B \subset {\cal X}_1^{m+1}\times{\cal X}_1^{2k+1}\times {\cal X}_{1}^{m+1}=
{\cal X}_1^{2m+2k+3}.$$ Clearly, the barriers from
$B$ need not be, and indeed, are not separated.  It is, however, easy to extend them to
separated ones. Indeed, let $q_0\ne 1$ be such that
$p_{q_0\,1}>0$ and redefine $B\stackrel{\mathrm{def}}{=}{\cal X}_{q_0}\times B$. Evidently,
any shift of any $y\in B$ by $w$ ($1\le w\le r$) positions to the right makes it impossible
for $y_1$ to be simultaneously in $\mathcal{X}_{q_0}$ and in $\mathcal{X}_1$ (since the latter sets are
disjoint, \S\ref{subsec:Xl}).
\end{proof}
\section{Conclusion}\label{sec:end}
As discussed in \S\ref{sec:intro} and \S\ref{sec:prework} in particular,  
the proper infinite alignments (\S\ref{sec:process}) allow us to define 
the decoding process $V$ which is regenerative and can further be 
stationarized to become ergodic \cite{AVT4}. This in turn allows us
to study the distribution and asymptotic properties not only 
of the Viterbi process $V$ but also of the joint process $(X,V)$.
In particular, this reveals how different 
these properties are from the properties of  the 
underlying chain $Y$ and HMM $(X,Y)$, respectively.
More specifically, since the process $V$ (resp. $(X,V)$)  
can deviate from the process $Y$ (resp. $(X,Y)$) significantly,
using  the Viterbi alignments $v_{1:n}$ as  estimates 
for  the hidden paths $Y_{1:n}$ might lead to incorrect conclusions
not only for finite $n$ (as generally appreciated) but also in the
limit as $n\to\infty$ \cite{AVT4}.  

This certainly does not mean  that one should not make 
inference based on $V$ but simply suggests that the aforementioned differences
may need to be taken into account. One example of how these asymptotic differences can be 
successfully accounted for  is the adjusted Viterbi training for 
HMM parameter estimation \cite{AVT1, AVT3, AVT4}.  

If known --- possibly estimated ---  these differences might also be appreciated
when the Viterbi paths are used for prediction, or segmentation, of $Y$,
e.g. in speech segmentation or in segmentation of DNA sequences 
into coding and non-coding regions, or in detection of CpG islands
in DNA sequences \cite{BioHMM2}. 
Indeed, in  segmentation of DNA sequences, the underlying chain $Y$ 
has few, often two, states (e.g. coding and non-coding
regions, or CpG islands and non-CpG regions), 
the probabilities of transitions between the states are very low,
hence the true ($Y$) and predicted ($V$) hidden paths consist of  
long constant blocks. At the same time, it has been noted 
that the predicted constant blocks can be somewhat longer 
than what the chain parameters would suggest. With the help of
the infinite Viterbi process $V$ it is now clear that this
discrepancy is not simply due to the random fluctuations but 
is systematic, does not vanish asymptotically, and is a direct
consequence of that  the transition probabilities of $V$ do indeed often
underestimate the true ones.  Note that in these examples, 
unlike in the estimation of the HMM emission parameters, the overall performance
is directly linked to the accuracy of the transition probability
estimates. Thus, finding the differences between the processes 
$(X,Y)$ and $(X,V)$ in this case might help find better alignments.


%

%
%
%
\ifCLASSOPTIONcompsoc
  \section*{Acknowledgments}
\else
  \section*{Acknowledgment}
\fi
The first author has been supported by the Estonian Science Foundation Grant 7553.
The authors thank Eurandom (The Netherlands) for initiating and stimulating 
their research on hidden Markov models, of which this work has been an 
integral part. The authors also thank Dr. A. Caliebe for valuable discussions
and for emphasizing the significance of the topic of path estimation 
in HMMs.

\ifCLASSOPTIONcaptionsoff
  \newpage
\fi



\bibliographystyle{IEEEtran}
%


\begin{IEEEbiography}{J\"uri  Lember} was born in 1968 in  Tallinn, Estonia.
He received the diploma and M.Sc. degrees in mathematical statistics in 1992 and 1994, respectively, 
from the University of Tartu, Estonia. He received the Ph.D. degree in mathematics from  
the University of Tartu, Estonia, in 1999. 

He completed his compulsory military service in 1987--1989, and was a Postdoctoral 
Research Fellow in the Institute of Mathematical Statistics, 
University of Tartu, in 1999--2000. He held a Postdoctoral Research position in Eurandom, 
The Netherlands, in 2001--2003. Since 2003, he has been a Lecture and a Senior Researcher 
in the Institute of Mathematical Statistics, University of Tartu. 
His  scientific interests include probability theory, theoretical 
statistics, information theory, and speech recognition.

Dr. Lember has been a member of the Estonian 
statistical society as well as Estonian mathematical society since 2003. He has
been awarded Estonian Science foundation grants for periods of 2004--2007 and 2008--2011.
\end{IEEEbiography}

\begin{IEEEbiography}{Alexey Koloydenko}
received the B.S. degrees in physics
and mathematics (with information systems minor) in 1994 from the Voronezh University, Russian Federation
and Norwich University, USA, respectively. He received in 1996 the M.S.(tech.)
degree in physics and radio-electronics from the Voronezh University, Russian Federation, 
and the M.S. degree in mathematics and statistics
from the University of Massachusetts at Amherst, USA. He received the Ph.D. degree
in mathematics and statistics from  the University of
Massachusetts at Amherst, USA, in 2000.

He held Postdoctoral Research and Teaching positions with the Department of Mathematics and Statistics
of the University of Massachusetts at Amherst, Statistics and Computer Science Departments of the University
of Chicago,  and Eurandom, The Netherlands, in 2000, 2001--2002, and 2002--2005, respectively. Since 2005 he
has been a Lecturer in Statistics at the University of Nottingham, UK. His research interests
include statistical processing and analysis of images, diffusion weighted MRI, 
algebraic aspects of probability theory and statistics, and hidden Markov models.

Dr. Koloydenko has been a member of the 
Pattern Analysis, Statistical Modelling and Computational Learning
European network (PASCAL) since 2004. 
\end{IEEEbiography}

\end{document}